\documentclass[11pt]{article}
\input epsf
\usepackage{amsfonts}
\usepackage{amscd}
\usepackage{amsmath,amssymb,amsthm}
\usepackage{amstext}
\usepackage{amscd}
\usepackage{graphicx}
\usepackage[all]{xy}
\usepackage{pstricks,pst-node,pst-tree}

\theoremstyle{plain}
\newtheorem{lem}{Lemma}[section]
\newtheorem*{class}{Theorem A}
\newtheorem*{case}{Theorem B}
\newtheorem*{weak}{Theorem C}
\newtheorem*{main}{Theorem 1}
\newtheorem*{main2}{Theorem 2}
\newtheorem*{main3}{Theorem 3}
\newtheorem*{main4}{Corollary 4}

\newtheorem{coro}[lem]{Corollary}

\theoremstyle{definition}

\newtheorem{definition}[lem]{Definition}

\newtheorem{rem}[lem]{Remark}
\newtheorem{ex}[lem]{Example}

\renewcommand{\descriptionlabel}[1]%
       {\hspace{\labelsep}\textsf{#1}}

\newcommand{\noi} {\noindent}

\DeclareMathOperator{\Area}{Area}
\DeclareMathOperator{\Spin}{Spin}
\DeclareMathOperator{\RCSpin}{RCSpin}
\DeclareMathOperator{\Long}{Length}

\begin{document}
\title{Minimal area of  the spun trefoil knot on the canonical cubulation of $\mathbb{R}^4$\,\,}
\author{ Ana Baray, Juan Jos\'e Catal\'an, Gabriela Hinojosa, Rogelio Valdez}
\date{July 12, 2025}

\maketitle
\begin{abstract} 
\noi \noi We say that a  \emph{cubical 2-knot}  $K^{2}$ is an embedding of the 2-sphere in the 2-skeleton of the canonical cubulation
of $\mathbb{R}^4$; in particular, $K^{2}$ is the union of $m(K^{2})$ unit squares, hence  $m(K^{2})$ is its area. We define the minimal area of $K^{2}$ as the minimum over all the areas of cubical 2-knots isotopic to the given knot type. The minimal area of a cubical 2-knot is an invariant, and the following natural question arose: Given a knot type, what area is needed for
a cubical 2-knot in the canonical cubulation of  $\mathbb{R}^4$ to realise that type with minimal area? 
In this paper, we answer this question for the spun trefoil knot in the weakly minimal case.
\end{abstract}

\rightline{Dedicated to Professor Alberto Verjovsky}


\section{Introduction}\label{sec1}

The study of polygons on the canonical cubulation of $\mathbb{R}^3$ started in the early 1960s (see \cite{diao}). Many interesting questions related to
the topological behavior of such polygons arose, for instance: What
is the minimal length for a polygon on the canonical cubulation of $\mathbb{R}^3$ that is knotted? Given a knot type, what is the minimal length needed for
a polygon on the canonical cubulation of  $\mathbb{R}^3$ to realise that knot type? Y. Diao proved in \cite{diao} that the minimal length for a polygon on the cubic lattice of $\mathbb{R}^3$ to be knotted is 24. We are interested in these questions for cubical 2-knots.\\

\noindent M. Boege, G. Hinojosa, and A. Verjovsky proved in \cite{BHV} the following result.\\

\begin{class}\label{class}
Any smooth $n$-knot ${K}^n:{\mathbb S}^n\hookrightarrow{\mathbb R}^{n+2}$ can be deformed isotopically into the $n$-skeleton  
of the canonical cubulation of  ${\mathbb R}^{n+2}$.\\
\end{class}

\noindent This isotopic copy is called a \emph{cubical $n$-knot}. As a consequence, every smooth 1-knot or 2-knot
$\mathbb{S}^n\subset{\mathbb {R}}^{n+2}$ ($n=1,2$) is isotopic to a cubical knot of the corresponding dimension. \\

\noindent The study of 2-knots in $\mathbb{R}^4$ has been considered by various authors, for instance in 
\cite{CRS}, \cite{CS}, \cite{K}, \cite{Roseman} and \cite{SSW}.\\

\noindent There are two types of elementary ``cubulated moves'' for cubical 2-knots (\cite{DHV}): 
The first one, $(M1)$, is obtained by dividing each hypercube of the original cubulation of 
$\mathbb{R}^4$ into $m^4$ hypercubes for $m\in\mathbb{N}$, and the second one $(M2)$ consists of exchanging a connected set of squared faces homeomorphic to a 
disk $\mathbb{D}^2$ in a cube of the cubulation with the complementary faces in that cube.\\

\noindent D\'iaz, Hinojosa and Verjovsky proved in \cite{DHV} the following.\\

\begin{case}\label{case}
Given two cubical 2-knots  $K^{2}_1$ and $K^{2}_2$  in $\mathbb{R}^{4}$, 
they are isotopic if and only if $K^{2}_1$ is equivalent to $K^{2}_2$ by a finite sequence of cubulated moves; {\emph{i.e.}}, $K^{2}_1\sim  K^{2}_2 \iff K^{2}_1\overset{c}\sim K^{2}_2$.\\
\end{case}

\noindent This result is analogous to the Roseman moves of classical tame 2-knots for cubical 2-knots.\\

\begin{figure}[h] 
\centering
 \includegraphics[height=3.5cm]{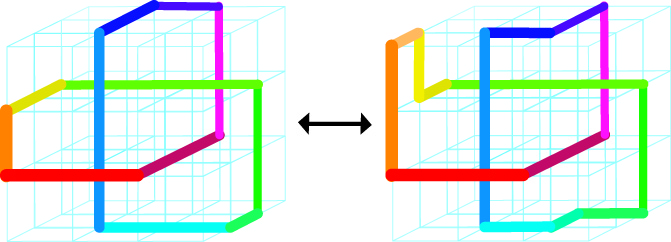}
\caption{\sl Two isotopic knots are equivalent via cubulated moves.} 
\label{M}
\end{figure} 

\noindent Notice that a cubical 2-knot $K^{2}$ is the union of $m(K^{2})$ unit squares, hence we have that $m(K^{2})$ is the \emph{area} of $K^{2}$. \\

\noindent Now, we define the \emph{minimal area of a cubical 2-knot} $K^{2}$, denoted by  $A(K^{2})$, as the minimum taken over all the areas of cubical 2-knots isotopic to $K^{2}$. By Theorem B, we have that $A(K^{2})$ is an invariant.\\

\noindent As a consequence, we can ask the following: Given a cubical knot type, what is the minimal area needed for
a cubical 2-knot on the canonical cubulation of  $\mathbb{R}^4$ to be a knot with the given knot type?  \\

\noindent We say that a cubical 2-knot $K^{2}$ is \emph{weakly minimal} if the area of $K^{2}$ is not reduced by any single cubulated move. In \cite{BH}, A. Baray and G. Hinojosa proved the next result.\\

\begin{weak}\label{weak}
There is a weakly minimal cubical 2-knot with area  256 representing the spun trefoil knot.\\
\end{weak}

\noindent Our goal in this paper is to prove that the weakly minimal area of the spun trefoil knot is 130, as stated below.\\

\begin{main}\label{main}
There is a weakly minimal cubical 2-knot with area  130 isotopic to the spun trefoil knot.\\
\end{main}

\noindent We will start studying the area of cubical spun 2-knots. Roughly speaking, a cubical spun 2-knot is constructed via the following process. Consider a cubical arc $Arc(K)$ in $\mathbb{R}^3_{+}$ whose end-points are in  $\mathbb{R}^{2}$. There exists a smooth arc $\overline{Arc(K)}$ isotopic to $Arc(K)$ such that $\overline{Arc(K)}$ is arbitrarily close to $Arc(K)$. Now, we spin the smooth arc $\overline{Arc(K)}$, to get  $\Spin(\overline{Arc(K)})\subset\mathbb{R}^4$.  Notice that $\Spin(\overline{Arc(K)})$ is a smooth 2-knot, hence  it is isotopic to a {\it cubical spun 2-knot}
$C\Spin(Arc(K))$ that can be obtained by gluing disks, annuli, and cylindrical surfaces properly, for more details see Section \ref{sec2.1} and \cite{BH}.
We can apply the same construction for  a $\lambda$-{\it subcubical}  $1$-arc $Arc(K)\subset \mathbb{R}^{3}_{+}$; \emph{i.e.,} 1-arc contained in the 1-skeleton of the subcubulation ${\cal{C}}_{\lambda}$,  for an integer $\lambda\geq1$ (see  (M1) - cubulated move); in particular, each vertex of $Arc(K)$ is contained in the 
lattice $(\frac{1}{\lambda}\mathbb{Z})^4$ and the length of each edge is $\frac{1}{\lambda}$. So, we obtain a $\lambda$-subcubical spun 2-knot 
$C\Spin(Arc(K))_{\lambda}$ contained in the 2-skeleton of the subcubulation ${\cal{C}}_{\lambda}$ of $\mathbb{R}^4$. \\

\begin{main2}\label{main2}
Let $C\Spin(Arc(K))_{\lambda}\subset {\cal{C}}_{\lambda}$ be a  $\lambda$-subcubical spun $2$-knot obtained from the $\lambda$-subcubical arc  $Arc(K)$ with vertices $v_{1}$, $v_{2}$, $\dots$, $v_{n}\in (\frac{1}{\lambda}\mathbb{Z})^4$, where $v_{i}=(x_{i},y_{i},z_{i},0)$ and $\frac{1}{\lambda}=||v_{i+1}-v_{i}||$. Then
 $$
 \Area(C\Spin(Arc(K))_{\lambda})=\frac{8}{\lambda}\sum_{i=2}^{n-1}z_{i}.
 $$
\end{main2}

\noindent As a consequence, we have that the minimal area of a cubical spun 2-knot $C\Spin(Arc(K))$ is upper bounded by 
$\displaystyle{8 \left(\sum_{i=2}^{n-1}z_{i} \right) -4n+6}$. \\

\noindent We denote the set of all cubical spun $2$-knots and reduced cubical spun $2$-knots by ${\mathcal{CS}}^2$. \\

\noindent The \emph{cubical spin minimal area} of a spun 2-knot $K^{2}$, denoted by  $A_{\mathcal{CS}}(K^{2})$, is defined as the minimal area of all cubical spun 2-knots and reduced cubical spun 2-knots isotopic to $K^{2}$; \emph{i.e.,}
$$
A_{\mathcal{CS}}(K^{2})=\min\{\Area(\tilde{K}^{2})\,:\,\tilde{K}^{2}\in {\mathcal{CS}}^2\,\,\,\mbox{and $\tilde{K}^2$ is isotopic to}\,\,K^{2}\}.
$$
Observe that the cubical minimal area of a spun 2-knot $K^2$ is greater than or equal to the minimal area of $K^2$; \emph{i.e.}, $A(K^{2})\leq A_{\mathcal{CS}}(K^{2})$.\\

\begin{main3}\label{main3}
The cubical spin minimal area of the spun trefoil knot $\Spin(T_{2,3})$ is 130, that is, 
$$
A_{\mathcal{CS}}(\Spin(T_{2,3}))=130.
$$
\end{main3}

\noindent Using these results, we prove our main theorem. As a consequence, the minimal area of the spun trefoil knot is smaller than or equal to 130.

\section{Preliminaries}

A cubulation of  $\mathbb{R}^{4}$ is a decomposition into a collection of right-angled $4$-dimensional
hypercubes $\{4,3,3\}$ called  $\textit{cells}$ such that any two are
either disjoint or meet in one common $k-$face of some dimension $k$. This provides  $\mathbb{R}^{4}$ with the structure of a cubical
complex whose category is similar to the simplicial category PL.  \\

\noindent The combinatorial structure of  any cubulation is the following: around each vertex, there are eight edges, 24 squares, 
32 cubes and 16 hypercubes. Around each edge, there are six squares, 12 cubes, and eight hypercubes. Around each square, there are four cubes and four hypercubes. Finally, around each cube, there are  
2 hypercubes (see  for instance \cite{DHV} and \cite{DHV1}).\\

\noindent The canonical cubulation $\mathcal{C} $ of $\mathbb{R}^{4}$ is its decomposition into hypercubes, which are the images of the unit hypercube
$$\{4,3,3\}=I^4=[0,1]^4=\{(x_{1},x_{2},x_{3},x_{4})\in \mathbb{R}^4 \,|\,0\leq x_{i}\leq 1\}$$ by translations by vectors with integer coefficients. 
Then, all vertices of $\mathcal{C} $ have integer coordinates. \\

\noindent Any regular hypercubic honeycomb $\{4,3,3,4\}$ or cubulation of $\mathbb{R}^{4}$ is obtained from the canonical cubulation by applying a conformal transformation  
to the canonical cubulation. Remember that a conformal transformation is of the form  
$x\mapsto{\lambda{A(x)}+a},$ where $\lambda\neq0,\,\,a\in \mathbb{R}^{4},\,\,\,A\in{SO(4)}$. \\

\begin{definition} The $k${\it-skeleton} of $\cal C$, denoted by $\mathcal{S}^k$, consists of the union of the $k$-skeletons of the hypercubes in  $\cal C$,
{\it i.e.,} the union of all cubes of dimension $k$ contained in the faces of the $4$-cubes in  $\cal C$.
We will call the 2-skeleton $\mathcal{S}^2$ of $\cal C$ the {\it canonical scaffolding} of $\mathbb{R}^{4}$.\\
\end{definition}

\noindent Notice that all the previous definitions can be extended naturally to $\mathbb{R}^{n+2}$ (compare \cite{BHV}).\\

\noindent A \emph{gridded surface} or {\it cubical surface} $S$ on $\mathcal{S}^2$ of $\cal C$ in $\mathbb{R}^{4}$ is a 
piecewise linear surface such that each linear piece is a unit square with its vertices in the $\mathbb{Z}^4$-lattice of $\mathbb{R}^{4}$ (see \cite{DHV1}).  
We note that our gridded surfaces are naturally length spaces \cite{BB}.

\subsection{Some facts of cubical 2-knots}

\noindent In classical knot theory, a subset $K$ of a space $X$ is a {\it knot} if $K$ is homeomorphic to a  $p$-dimensional sphere 
$\mathbb{S}^{p}$ embedded in either the Euclidean $n$-space $\mathbb{R}^{n}$ or the $n$-sphere $\mathbb{S}^{n}=\mathbb{R}^{n}\cup\{\infty\}$, where $p<n$. Two knots $K_1$, $K_2$ are 
{\it equivalent} or {\it isotopic} if there exists a homeomorphism $h:X\hookrightarrow X$ such that $h(K_1)=K_2$;
in other words $(X,K_1)\cong (X,K_2)$. However, a knot $K$ is sometimes defined to be an embedding
$K:\mathbb{S}^{p}\hookrightarrow\mathbb{S}^{n}\cong\mathbb{R}^{n} \cup \{\infty \}$ (see \cite{mazur}, \cite{rolfsen}).
We shall also find this convenient at times and will use the same symbol to
denote either the map $K$ or its image $K(\mathbb{S}^{p})$ in $\mathbb{S}^{n}$.\\

\begin{definition}
Let $K^2$ be a $2$-dimensional knot in $\mathbb{R}^{4}$. If $K^2$ is contained in the canonical scaffolding $\mathcal{S}^2$, we say that $K^{2}$ is a \emph{cubical 2-knot}. \\
\end{definition}



\begin{definition}
The \emph{minimal area} of  cubical 2-knot  $K^{2}$, denoted by  $A(K^{2})$, is defined as 
$$
A(K^{2})=\min\,\{\,\Area(T^{2})\,:\,T^{2}\,\,\,\mbox{is a cubical 2-knot isotopic to}\,\,K^{2}\,\}.
$$
\end{definition}

\noindent Observe that  \emph{minimal area} is an invariant for cubical 2-knots.\\

\begin{definition}
Given a cubical 2-knot $K^{2}$, we say that $K^{2}$ is \emph {minimal} if the area of $K^{2}$ is $A(K^{2})$.\\
\end{definition}

\begin{rem}
By Theorem A, we can extend the definition of minimal area  to any smooth 2-knot $K^{2}$, denoted by  $A(K^{2})$, as 
$$
A(K^{2})=\min\{\Area(\tilde{K}^{2})\,:\,\tilde{K}^{2}\,\,\,\mbox{is a cubical knot isotopic to}\,\,K^{2}\}.
$$
Notice that by Theorem B, $A(K^{2})$ is also an invariant.\\
\end{rem}

\begin{ex}
The minimal area of the unknotted sphere is 6.
\end{ex}

\subsection{Cubulated moves}

In cubical 1-knots or lattice knots, there are two types of elementary ``cubulated moves''. The first one (M1) is obtained by dividing each cube of the original cubulation of $\mathbb{R}^3$ into $m^3$ cubes, which means that each edge of the knot is subdivided into $m$ equal segments. The second one (M2) consists of exchanging a connected set of edges in a face of the cubulation with the complementary edges in that face (see \cite{HVV}). These cubulated moves can be naturally extended to cubical 2-knots; for more details, see \cite{DHV}. \\

\begin{definition} The following are the allowed {\emph{cubulated moves}}:
\begin{description}
\item[{\bf{M1}}] {\emph{Subdivision:}} Given an integer $m\geq 1$, consider the subcubulation  ${\cal{C}}_m$ of $\cal{C}$ by subdividing each $k$-dimensional cell of $\cal{C}$ in 
$m^k$ congruent $k$-cells in ${\cal{C}}_m$, in particular each hypercube in $\cal{C}$ is subdivided in $m^4$ congruent hypercubes in ${\cal{C}}_m$. Moreover, as a cubical complex,
each $k$-dimensional face of the 2-knot $K^2$ is subdivided into $m^k$ congruent 
$k$-faces. Since ${\cal{C}}\subset {\cal{C}}_{m}$, then
$K^2$ is contained in the scaffolding ${\cal{S}}^2_m$ (the $2$-skeleton) of ${\cal{C}}_m$. \\
 
\item[{\bf{M2}}] {\emph{Face Boundary Moves:}} Suppose that $K^2$ is contained in some subcubulation ${\cal{C}}_m$ of  
the canonical cubulation ${\cal{C}}$ of $\mathbb{R}^{4}$. 
Let $Q^4\in {\cal{C}}_m$ be a $4$-cube such that
$A^2=K^2\cap Q^4$ contains a $2$-face. We can assume, up to applying the elementary $(M1)$-move if necessary, that $A^2$ consists of either one, 
two, or three
squares such that it is a connected surface and it is contained in the boundary of a 3-cube
$F^3\subset Q^4$; in other words, $A^2$ is a cubical disk contained in the boundary of
$F^3$. The boundary $\partial F^3$ is divided by $\partial A^2$ into two cubulated surfaces, one of which is $A^2$, and we denote the other by $B^2$. Observe
that both cubulated surfaces share a common circle boundary. The face boundary move consists in replacing $A^2$ by $B^2$ (see Figure \ref{M1}). There are four types of face
boundary moves depending on the number of 2-faces in each $A^2$ and $B^2$. If $A^2$ has $p$ squares then $B^2$ has $6-p$ squares. 

\begin{figure}[h] 
\centering
 \includegraphics[height=5cm]{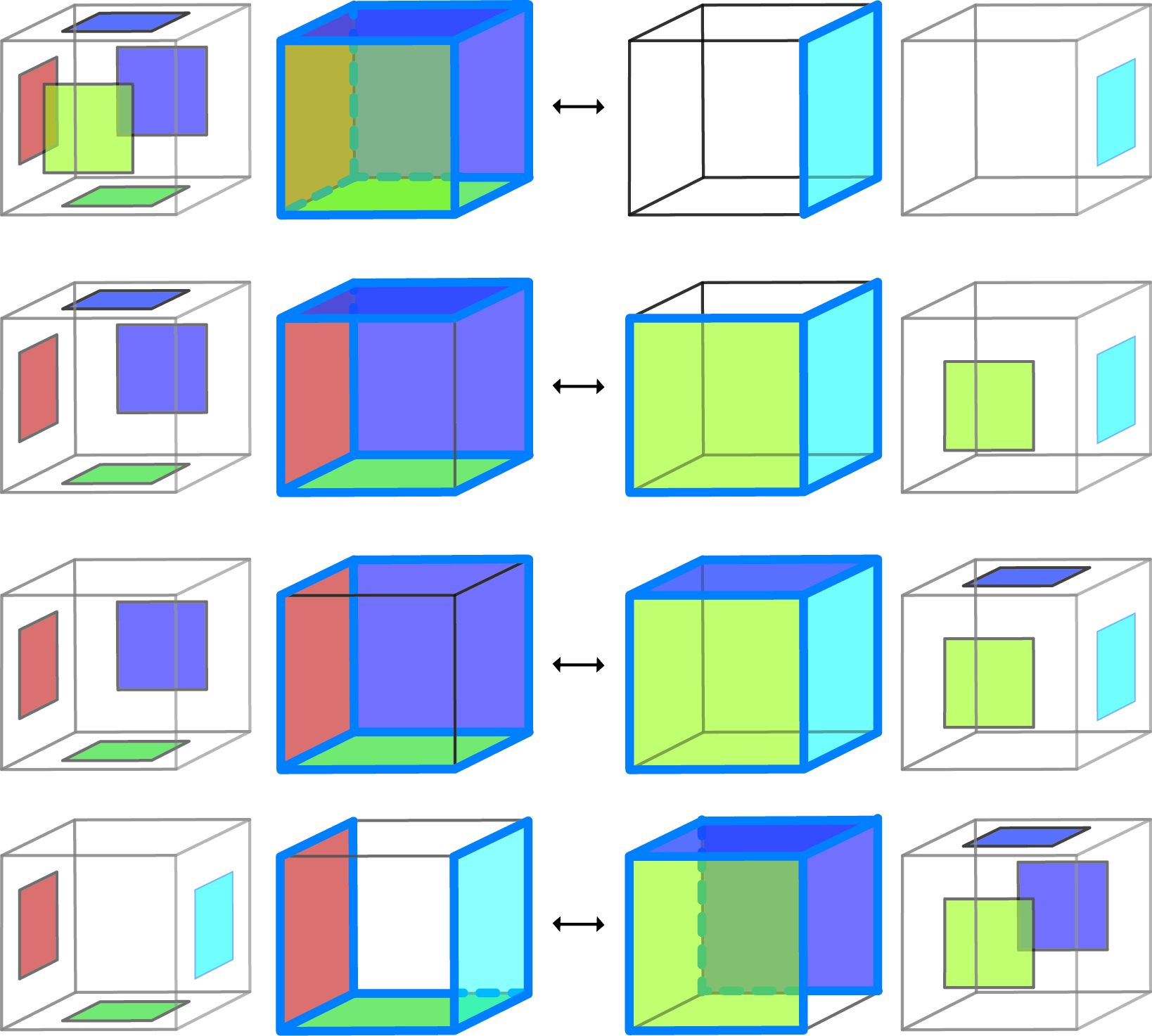}
\caption{\sl The four types of face boundary moves.} 
\label{M1}
\end{figure} 
\end{description}
\end{definition}

\begin{rem}\label{equivalencia}
Notice that  the $(M2)$-move can be extended to an ambient isotopy of $\mathbb{R}^4$.\\
\end{rem}

\begin{definition}
Consider two cubical $2$-knots $K^2_1$ and $K^2_2$ in $\mathbb{R}^{4}$. We say that
$K^2_1$ is {\emph{equivalent}} to $K^2_2$ {\emph{by cubulated moves}}, denoted by 
$K^2_1\overset{c}\sim K^2_2$, if we can transform $K^2_1$ into $K^2_2$ by a finite number
of cubulated moves. \\
\end{definition}

\noindent The following result was proved in \cite{DHV}. \\

\begin{case}\label{case}
Given two cubical 2-knots  $K^2_1$ and $K^2_2$  in $\mathbb{R}^{4}$, 
they are isotopic if and only if $K^2_1$ is equivalent to $K^2_2$ by a finite sequence of cubulated moves; {\emph{i.e.}}, $K^2_1\sim  K^2_2 \iff K^2_1\overset{c}\sim K^2_2$.\\
\end{case}

\begin{definition}
Given a cubical 2-knot $K^{2}$, we say that $K^{2}$  is \emph{weakly minimal} if the area of $K^{2}$ can not  be reduced by any single cubulated move.\\
\end{definition}

\noindent Clearly, if a 2-knot is minimal, then it is weakly minimal.

\section{The area of cubical spun 2-knots}\label{sec2}

The purpose of this section is to study the area of the spun trefoil 2-knot. We will start constructing cubical spun 2-knots in an analogous way to that
used in  \cite{BH}.\\

\noindent Consider in $\mathbb{R}^{4}$ the half-space
$$
\mathbb{R}^{3}_{+}=\{(x_{1},x_{2},x_{3},0):x_{3}\geq 0\},
$$
whose boundary is the plane $\mathbb{R}^{2}=\{(x_{1},x_{2},0,0)\}.$\\

\noindent We can spin each point $x=(x_{1},x_{2},x_{3},0)$ of $\mathbb{R}^{3}_{+}$
with respect to $\mathbb{R}^{2}$ according to the formula
$$
R_{\theta}(x)=(x_{1},x_{2},x_{3}\cos (2\pi \theta),x_{3}\sin(2\pi \theta)).
$$
\noindent We define $\Spin(X)$ of a set $X\subset\mathbb{R}^{3}_{+}$, as
$$
\Spin(X)=\{R_{\theta}(x):x\in X, 0\leq\theta\leq 1\}.
$$
\noindent Given a tame arc $A$ in $\mathbb{R}^{3}_{+}$ such that its end-points are in  $\mathbb{R}^{2}$ and its interior is  in 
$\mathbb{R}^{3}_{+}-\mathbb{R}^{2}$, we can obtain a 2-knot 
$\Spin(A)$ in $\mathbb{R}^{4}$ called it  {\it spun knot}.\\

\noindent Notice that  we can consider $A$ as the image of an embedding of the unit interval,
$A:I\rightarrow\mathbb{R}^{3}_{+}$, such that $A(0)\neq A(1)\in\mathbb{R}^{2}$ and $A(0,1)\subset \mathbb{R}^{3}_{+}-\mathbb{R}^{2}$, which will
be denoted by the same letter. Then we will say that an arc $A\subset\mathbb{R}^{3}_{+}$
is a {\it spinnable arc} if it is a proper arc, i.e., it is smooth at every interior point and it has contact of infinite order at each end-point concerning the corresponding normal vector.\\

\subsection{Construction}\label{sec2.1}

\noindent Consider a cubical $1$-knot  $K\subset \mathbb{R}^{3}_{+}$, then each vertex of $K$ belongs to the lattice $\mathbb{Z}^4$, and each edge is parallel to some coordinate axis and has length one. Thus, we obtain a cubical arc from it, $Arc(K)$ with vertices  $v_{1}$, $v_{2}$, $\dots$, $v_{n}$ in 
$\mathbb{R}_{+}^3$, and edges  $\overline{v_{1}v_{2}}$, $\overline{v_{2}v_{3}}$, $\dots$, $\overline{v_{n-1}v_{n}}$, where $v_{1}$ and $v_{n}$ are in $\partial\mathbb{R}_{+}^{3}=\mathbb{R}^{2}$, and $v_{2}$, $v_{3}$, $\dots$, $v_{n-1}$ are in $\mathbb{R}_{+}^{3}-\mathbb{R}^{2}$. There exists a 
smooth arc $\overline{Arc(K)}$ isotopic to $Arc(K)$ such that $\overline{Arc(K)}$ is 
arbitrarily close to $Arc(K)$. This is because
we can round the corners at the vertices of $Arc(K)$ in arbitrarily small neighborhoods of them (see \cite{douady}).
Now, we will apply the map $\psi: \mathbb{R}^3_{+}\times [0,1]\rightarrow \mathbb{R}^4$ given by $(x,y,z,0, \theta)\mapsto (x,y,z\cos (2\pi\theta),z\sin (2\pi\theta))$ to obtain the spin of $\overline{Arc(K)}$, $\Spin(\overline{Arc(K)})\subset\mathbb{R}^4$. \\

\begin{figure}[h] 
\begin{center}
\includegraphics[width=12cm, height=4cm ]{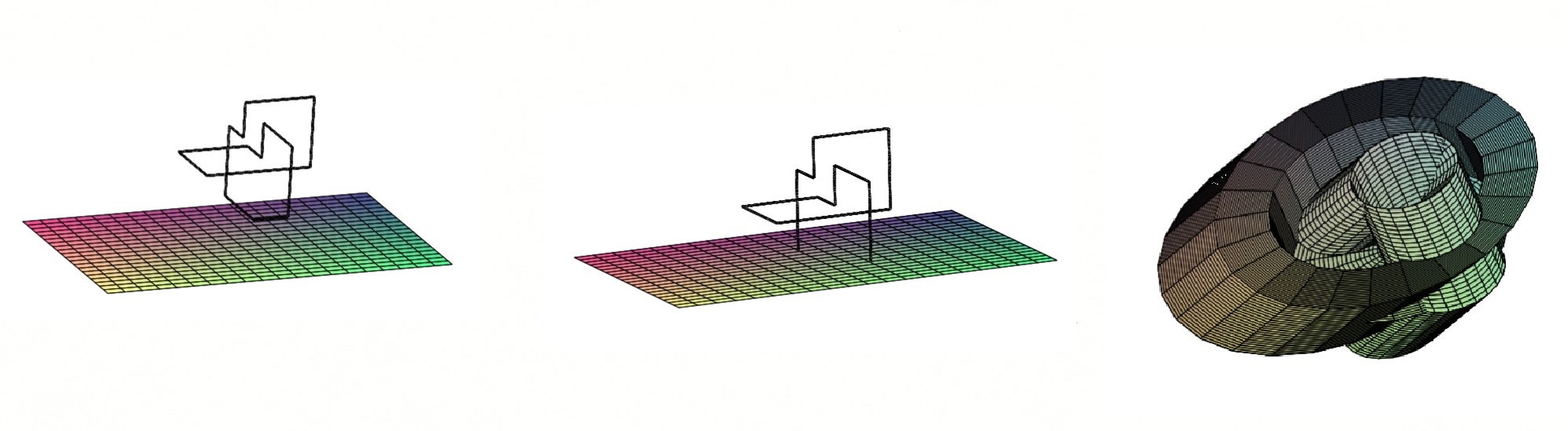}
\caption{\sl A decomposition of $\Spin\big(\overline{Arc(K)}\big)$ by pieces.}
\label{spintrefoil}
\end{center}
\end{figure}

\noindent Observe that $\Spin(\overline{Arc(K)})$ is a smooth 2-knot, hence by \cite{BHV}, it is isotopic to a cubical 2-knot $C\Spin(Arc(K))$. Since  
$\overline{Arc(K)}$ and $Arc(K)$ are arbitrarily closed, and $\Spin(Arc(K))=\cup_{i=1}^{n-1}\Spin(\overline{v_{i}v_{i+1}})$, we can think  $\Spin(\overline{Arc(K)})$ as a  $2$-sphere embedded in $\mathbb{R}^{4}$ obtained by gluing disks, annuli, and cylindrical surfaces properly.  We illustrate this decomposition in Figure \ref{spintrefoil}.\\

\noindent Next, we will construct the cubical 2-knot $C\Spin(Arc(K))$ from $\Spin(\overline{Arc(K)})$. 
The main idea is to cubulate each piece $\Spin(\overline{v_{i}v_{i+1}})$ for any $i$, in such a way that each step is compatible with the previous one. This procedure is analogous to that applied in  \cite{BH}.\\

\noindent Consider the canonical vectors in  $\mathbb{R}^{4}$: $e_{\pm1}=(\pm1,0,0,0)$, $e_{\pm2}=(0,\pm1,0,0)$, $e_{\pm3}=(0,0,\pm1,0)$, $e_{\pm4}=(0,0,0,\pm1)$. Then $v_{i+1}-v_{i}=e_{k_{i}}$, where $k_{i} \in \{ \pm1, \pm2, \pm3 \}$. For each vertex $v_{i}=(x_{i},y_{i},z_{i},0)\in Arc(K)$, let $f_{i}:[0,1] \rightarrow \mathbb{R}^{4}$ be the function defined as $f_{i}(s)=v_{i}+s(v_{i+1}-v_{i})$. Then, the image of the closed interval  $[0,1]$ under $f_{i}$ is the edge 
$\overline{v_{i}v_{i+1}}$, hence
\begin{eqnarray*} 
\Spin(\overline{v_{i}v_{i+1}}) &=& \{R_{\theta}(x) \,|\, x\in \overline{v_{i}v_{i+1}}, \, 0\leq \theta \leq 1  \}  \\
                               &=& \{R_{\theta}(f_{i}(s)) \,|\, 0\leq s \leq 1, \, 0\leq \theta \leq 1  \} \\
															 &=& \{R_{\theta}(f_{i}(1-s)) \,|\, 0\leq s \leq 1, \, 0\leq \theta \leq 1  \}.
\end{eqnarray*}

\noindent In order to describe  $\Spin(\overline{v_{i}v_{i+1}})$, we will analyze each case. \\

\noindent {\it Case 1.} The edge $v_{i+1}-v_{i}=e_{3}$, so $f_{i}(s) = v_{i}+se_{3}$, and
$$
\Spin(\overline{v_{i}v_{i+1}}) = \left \{\big(x_{i},y_{i},(z_{i}+s)\cos (2\pi\theta), (z_{i}+s)\sin(2\pi\theta) \big)  
\left |
\begin{array}{l}
0\leq s \leq 1, \\
0\leq \theta \leq 1 
\end{array}
\right.																
\right\}.
$$

\begin{figure}[h] 
\begin{center}
\includegraphics[width=10cm, height=3.5cm ]{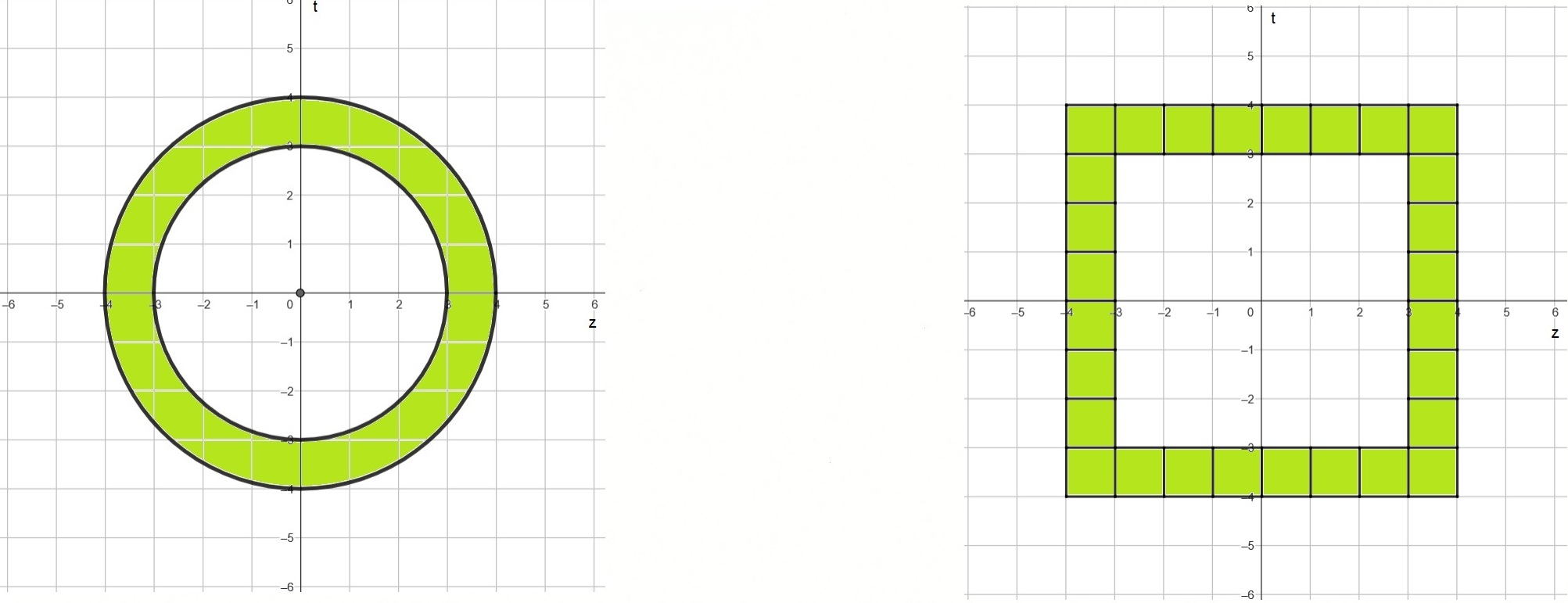}
\caption{\sl An annulus can be deformed into a square annulus.}
\label{annulus}
\end{center}
\end{figure}

\noindent If $z_{i}=0$, we get a disk centered at $(x_{i},y_{i},0,0)$ with radius $1$. If $z_{i}>0$, we obtain an annulus centered at $(x_{i},y_{i},0,0)$, with smallest radius $z_{i}$ and greatest radius equal to $z_{i}+1$.  Then $\Spin(\overline{v_{i}v_{i+1}})$ can be deformed into a square annulus as it is shown in Figure \ref{annulus}.\\

\noindent {\it Case 2.} The edge $v_{i+1}-v_{i}=e_{-3}$, so $f_{i}(1-s) = v_{i+1}+se_{3}$, and
$$
\Spin(\overline{v_{i}v_{i+1}}) =\,\left \{\big(x_{i+1},y_{i+1},(z_{i+1}+s)\cos (2\pi\theta), (z_{i+1}+s)\sin (2\pi\theta) \big) 
\left |
\begin{array}{l}
0\leq s \leq 1, \\
0\leq \theta \leq 1 
\end{array}
\right.																
\right\}. 
$$
\noindent If $z_{i+1}=0$, we get a disk centered at $(x_{i+1},y_{i+1},0,0)$ with radius $1$, and if  $z_{i+1}>0$, we obtain an annulus centered at  
$(x_{i+1},y_{i+1},0,0)$ with smallest radius equal to $z_{i+1}$ and greatest radius equal to $z_{i+1}+1$.  As above,  $\Spin(\overline{v_{i}v_{i+1}})$ can be deformed into a square annulus (see Figure \ref{annulus}).\\

\begin{figure}[h] 
\begin{center}
\includegraphics[width=10cm, height=3.5cm ]{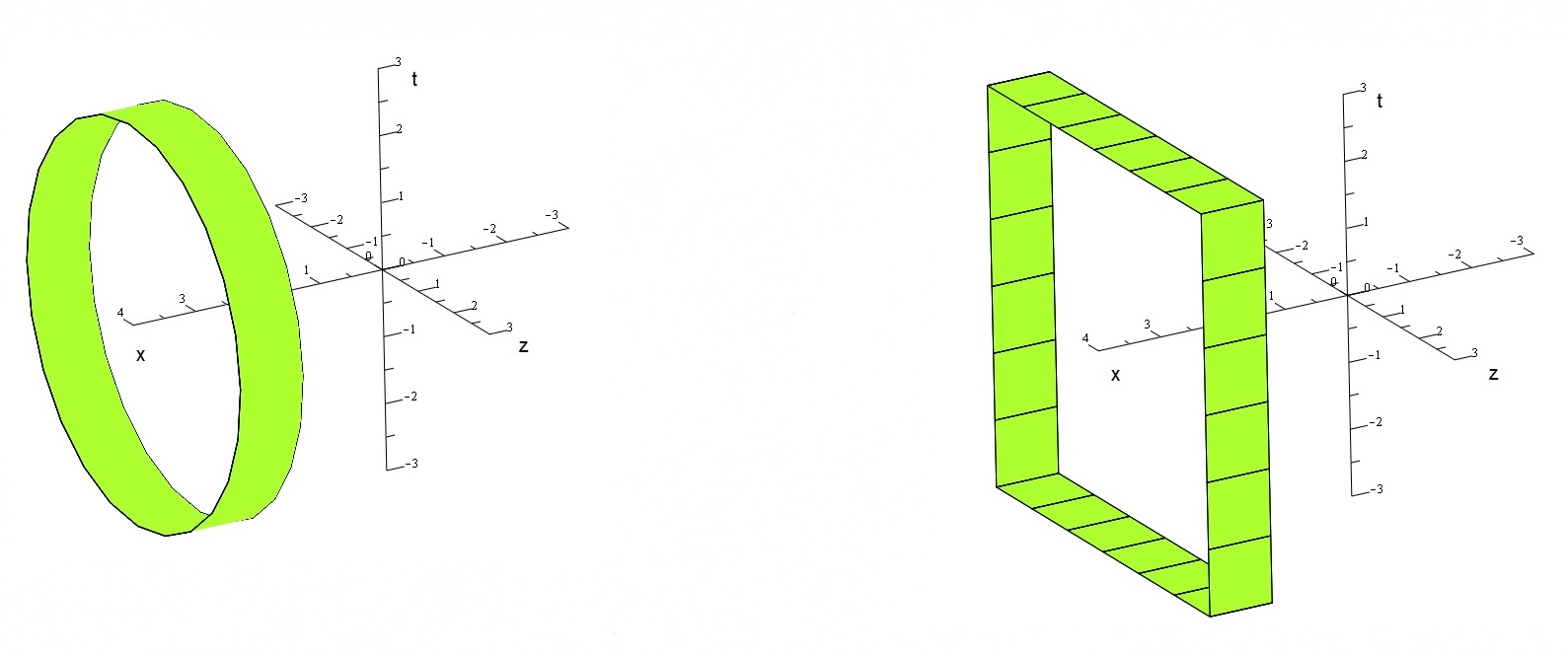}\\
\caption{\sl A cylinder can be deformed into a square cylinder.}
\label{cylinder}
\end{center}
\end{figure}

\noindent {\it Case 3.}  The edge $v_{i+1}-v_{i}=e_{1}$. Consider the function $f_{i}(s) = v_{i}+se_{1}$, so
$$
\Spin(\overline{v_{i}v_{i+1}}) = \{(x_{i}+s,y_{i},z_{i}\cos (2\pi\theta), z_{i}\sin (2\pi\theta) ) \,|\, 0\leq s \leq 1, \, 0\leq \theta \leq 1  \}.
$$
\noindent Then, we obtain a cylinder with a radius of $z_{i}$ and a height of 1 in a direction parallel to the $x$-axis. Its bottom boundary component is a circle centered at $(x_{i},y_{i},0,0)$, and its top boundary component is a circle centered at $(x_{i}+1,y_{i},0,0)$. Hence $\Spin(\overline{v_{i}v_{i+1}})$ can be deformed into a square cylinder without the bottom and top covers, as it is shown in Figure \ref{cylinder}.\\

\noindent {\it Case 4.} The edge $v_{i+1}-v_{i}=e_{-1}$. We have that $f_{i}(1-s) = v_{i+1}+se_{1}$, then 
$$
\Spin(\overline{v_{i}v_{i+1}}) = \left \{(x_{i+1}+s,y_{i+1},z_{i+1}\cos (2\pi\theta), z_{i+1}\sin (2\pi\theta) )  
\left |
\begin{array}{l}
0\leq s \leq 1, \\
0\leq \theta \leq 1 
\end{array}
\right.																
\right\}.
$$

\noindent Hence, we obtain a cylinder with radius  $z_{i+1}$ and height equal to 1 in a direction parallel to the coordinate axis $x$. 
 As above,  $\Spin(\overline{v_{i}v_{i+1}})$ can be deformed into a square cylinder without the bottom and top covers (see Figure \ref{cylinder}).\\

\noindent {\it Case 5.} The edge $v_{i+1}-v_{i}=e_{2}$. So $f_{i}(s) = v_{i}+se_{2}$, and
$$
\Spin(\overline{v_{i}v_{i+1}}) = \{(x_{i},y_{i}+s,z_{i}\cos (2\pi\theta), z_{i}\sin (2\pi\theta) ) \,|\, 0\leq s \leq 1, \, 0\leq \theta \leq 1  \}.
$$

\begin{figure}[h] 
\begin{center}
\includegraphics[width=10cm, height=3.5cm ]{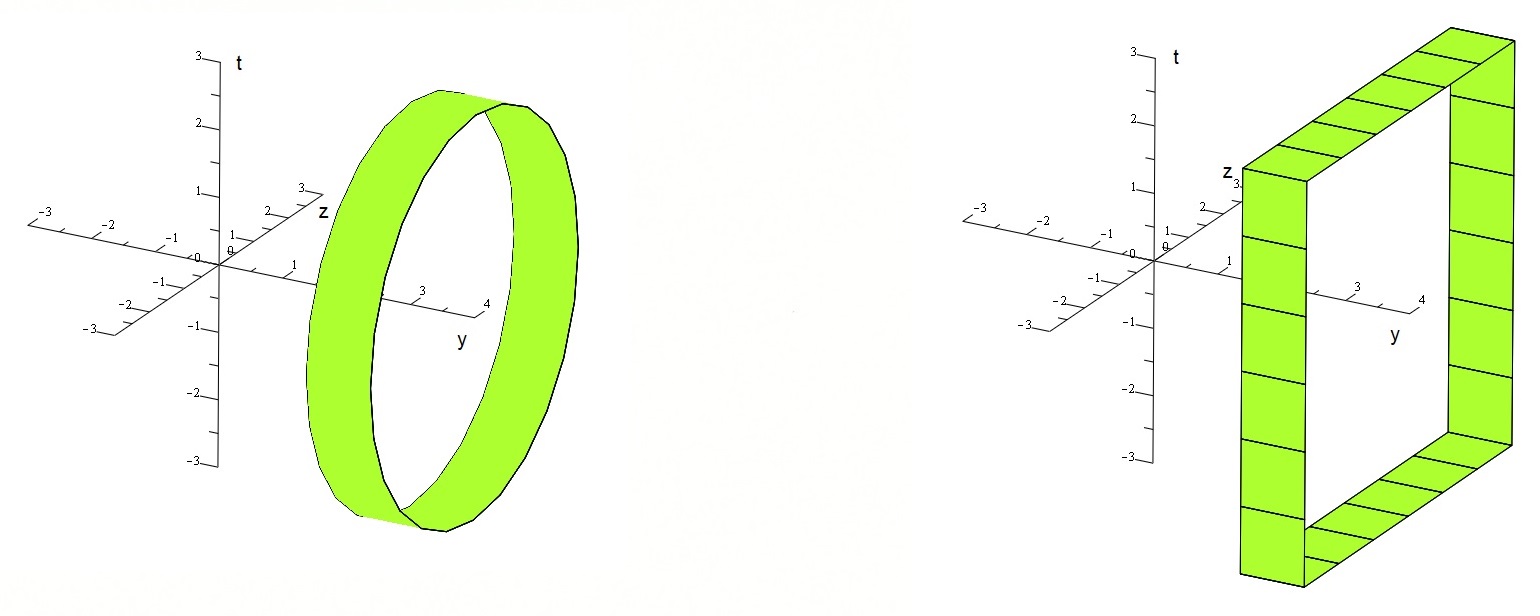}
\caption{\sl A cylinder can be deformed into a square cylinder.}
\label{rectangular}
\end{center}
\end{figure}

\noindent As above, we obtain a cylinder with radius  $z_{i}$ and a height of 1 in a direction parallel to the $y$-coordinate axis. 
Its bottom boundary component is a circle centered at $(x_{i},y_{i},0,0)$, and its top boundary component is a circle centered at  $(x_{i},y_{i}+1,0,0)$. Again,
 $\Spin(\overline{v_{i}v_{i+1}})$ can be deformed into a square cylinder without the bottom and top covers (see Figure \ref{rectangular}).\\

\noindent {\it Case 6.} The edge $v_{i+1}-v_{i}=e_{-2}$. As in the previous case $f_{i}(1-s) = v_{i+1}+se_{2}$, this implies that 
$$
\Spin(\overline{v_{i}v_{i+1}}) = \left \{(x_{i+1},y_{i+1}+s,z_{i+1}\cos (2\pi\theta), z_{i+1}\sin (2\pi\theta))  
\left |
\begin{array}{l}
0\leq s \leq 1, \\
0\leq \theta \leq 1 
\end{array}
\right.																
\right\}.
$$

\noindent Again, we obtain a cylinder with a radius of $z_{i+1}$ and a height of 1 in a direction parallel to the $ y$-axis.  Again,
 $\Spin(\overline{v_{i}v_{i+1}})$ can be deformed into a square cylinder without the bottom and top covers (see Figure \ref{rectangular}).\\
 
 \noindent Observe that each step of our construction is compatible with the previous one. \\

\begin{definition}
We say that a $2$-knot $K^{2}$ is a \emph{cubical spun $2$-knot}, if it can be obtained from a cubical arc via the process described above. \\
\end{definition}


\begin{rem}\label{subcubical}
Suppose that we start with a $\lambda$-{\it subcubical}  $1$-knot  $K\subset \mathbb{R}^{3}_{+}$; \emph{i.e.,} 1-knot contained in the 1-skeleton of the subcubulation ${\cal{C}}_{\lambda}$,  for some integer $\lambda\geq 1$ (see  (M1)- cubulated move in Section 2.2). In particular, each vertex of $K$ is 
contained in the lattice $(\frac{1}{\lambda}\mathbb{Z})^4$ and the length of each edge is $\frac{1}{\lambda}$. Using the previous procedure, we can construct a {\it subcubical spun 2-knot} $\Spin(Arc(K))_{\lambda}\subset {\cal{C}}_{\lambda}$ spinning the arc $Arc(K)$. More precisely, let $Arc(K)$ be the arc obtained from $K$, with vertices  $v_{1}$, $v_{2}$, $\ldots$, $v_{n}$ in $\mathbb{R}_{+}^{3}$, and edges
$\overline{v_{1}v_{2}}$, $\overline{v_{2}v_{3}}$, $\ldots$, $\overline{v_{n-1}v_{n}}$, such that $v_{1}$ and $v_{n}$ are in 
$\partial\mathbb{R}_{+}^{3}=\mathbb{R}^{2}$, and $v_{2}$, $v_{3}$, $\ldots$, $v_{n-1}$ are in the interior of  $\mathbb{R}_{+}^{3}$. As above, we have that 
$\Spin\big(Arc(K)\big)_{\lambda}=\cup_{i=1}^{n-1}\Spin(\overline{v_{i}v_{i+1}})$ 
is a  $2$-sphere embedded in $\mathbb{R}^{4}$ obtained by gluing disks, annuli, and cylindrical surfaces in a proper way, where $v_{i+1}-v_{i}=\frac{1}{\lambda} e_{k_{i}}$, for $k_{i} \in \{ \pm1, \pm2, \pm3 \}$, $\frac{1}{\lambda}=||v_{i+1}-v_{i}||$. Then, using the same procedure,  we construct a subcubical spun 2-knot $C\Spin(Arc(K))_{\lambda}$ from $\Spin(\overline{Arc(K)})_{\lambda}$ cubulating each piece. Still, this 2-knot is contained in the 2-skeleton of the subcubulation ${\cal{C}}_{\lambda}$ of $\mathbb{R}^4$. \\
\end{rem}

\begin{definition}
We say that a $2$-knot  $K^2_{\lambda}$ is a \emph{$\lambda$-subcubical spun $2$-knot} if it is contained in the 2-skeleton of the subcubulation ${\cal{C}}_{\lambda}$, and it can be obtain from a $\lambda$-subcubical arc $Arc(K)\subset {\cal{C}}_{\lambda}$, via the  process described above. \\
\end{definition}

\noindent We denote the set of all $\lambda$-subcubical spun $2$-knots contained in the subcubulation ${\cal{C}}_{\lambda}$ by ${\mathcal{CS}}^2_{\lambda}$, 
for an integer $\lambda\geq 1$.
\subsection{Area of subcubical and cubical spun 2-knots}

In this section, we will obtain a formula to compute the area of a cubical spun 2-knot $C\Spin(Arc(K))$, which will depend only on the arc $Arc(K)$. We will start considering $\lambda$-subcubical spun 2-knots.\\

\noindent Let $C\Spin(Arc(K))_{\lambda}\subset {\cal{C}}_{\lambda}$ be the subcubical spun $2$-knot obtained from the subcubical arc $Arc(K)$ (see Remark \ref{subcubical}). Let $v_{1}$, $v_{2}$, $\dots$, $v_{n}$ be the corresponding vertices of  $Arc(K)$, then the subcubulated surface generated by the edge 
$\overline{v_{i}v_{i+1}}$ will be denoted by $S_{i}$. Notice that the intersection of $S_{i}$ and $S_{j}$, $i\neq j$, is either the empty set or a boundary component, hence
$$
\Area(C\Spin(Arc(K))_{\lambda})=\sum_{i=1}^{n-1}\Area(S_{i}). 
$$
\noindent Next, we will compute the area of  $S_{i}$ in terms of the corresponding vertices $v_{i}$, $v_{i+1}$, where $v_{i}=(x_{i},y_{i},z_{i},0)$.\\

\noindent Suppose that $v_{i+1}-v_{i}=\frac{1}{\lambda} e_{k_{i}}$, $k_{i}=\pm 3$, then $||v_{i+1}-v_{i}||=|z_{i+1}-z_{i}|=\frac{1}{\lambda}$. 
If $k_{i}=3$, then the surface $S_{i}$ is a squared annulus (see Figure \ref{annulus}),  the length of the side of the smaller square is equal to $2z_{i}$, and the length of the side of the bigger square is equal to  $2z_{i+1}$. Notice that if $z_{i}=0$, then $S_{i}$ is a square. Now, for $k_{i}=-3$ we have again that  $S_{i}$ is a squared annulus such that the length of the side of the smaller square is equal to $2z_{i+1}$ and the length of the side of the bigger square is $2z_{i}$; and again if $z_{i+1}=0$, it follows that $S_{i}$ is a square. In both cases, we have that.
$$
\begin{array}{lllll}
\Area(S_{i}) &=& |(2z_{i})^{2}-(2z_{i+1})^{2}| \\
                  &=& 4|z_{i}^{2}-z_{i+1}^{2}|  \\
                  &=& 4(z_{i}+z_{i+1})|z_{i}-z_{i+1}|\\
                  &=& \frac{4}{\lambda} (z_{i}+z_{i+1}).
\end{array}
$$
\noindent
Suppose that $v_{i+1}-v_{i}=\frac{1}{\lambda} e_{k_{i}}$, $k_{i}=\pm 1$. In both cases, $||v_{i+1}-v_{i}||=|x_{i+1}-x_{i}|=\frac{1}{\lambda}$ and $z_{i}=z_{i+1}$. This implies that the cubulated surface $S_{i}$ is a square cylinder without the bottom and top covers  (see Figure \ref{cylinder}), such that the base length of each rectangle is $2z_{i}$. The length  of its height is $|x_{i+1}-x_{i}|$. Hence, 
$$
\begin{array}{lllll}
\Area(S_{i}) &=& 4(2z_{i})|x_{i+1}-x_{i}| \\
                  &=& 4(z_{i}+z_{i+1})|x_{i+1}-x_{i}| \\
                  &=& \frac{4}{\lambda} (z_{i}+z_{i+1}).
\end{array}
$$
\noindent
Finally, consider that $v_{i+1}-v_{i}=\frac{1}{\lambda} e_{k_{i}}$, $k_{i}=\pm 2$. In any of these cases $||v_{i+1}-v_{i}||=|y_{i+1}-y_{i}|=\frac{1}{\lambda}$, and as above $S_{i}$ is a 
square cylinder without the bottom and top covers  (see Figure \ref{rectangular}). As above
$$
\begin{array}{lllll} 
\Area(S_{i}) &=& 4(z_{i}+z_{i+1})|y_{i+1}-y_{i}| \\
                  &=& \frac{4}{\lambda} (z_{i}+z_{i+1}).
\end{array}
$$
\noindent
\begin{main2}\label{main2}
Let $C\Spin(Arc(K))_{\lambda}\subset {\cal{C}}_{\lambda}$ be a  $\lambda$-subcubical spun $2$-knot obtained from the $\lambda$-subcubical arc  $Arc(K)$ with vertices $v_{1}$, $v_{2}$, $\dots$, $v_{n}\in (\frac{1}{\lambda}\mathbb{Z})^4$, where $v_{i}=(x_{i},y_{i},z_{i},0)$ and $\frac{1}{\lambda}=||v_{i+1}-v_{i}||$. Then
 $$
 \Area(C\Spin(Arc(K))_{\lambda})=\frac{8}{\lambda}\sum_{i=2}^{n-1}z_{i}.
 $$
\end{main2}
\noindent {\it Proof.} By the previous discussion, we have that 
\begin{eqnarray*} 
\Area(C\Spin(Arc(K))_{\lambda}) &=& \frac{4}{\lambda}\sum_{i=1}^{n-1}(z_{i}+z_{i+1}) \\
             &=& \frac{4}{\lambda}\sum_{i=1}^{n-1}z_{i}+\frac{4}{\lambda} \sum_{i=1}^{n-1}z_{i+1}  \\
             &=& \frac{4}{\lambda}\sum_{i=1}^{n-1}z_{i}+\frac{4}{\lambda} \sum_{i=2}^{n}z_{i}  \\
             &=& \frac{4}{\lambda} z_{1}+\frac{8}{\lambda}\sum_{i=2}^{n-1}z_{i}+\frac{4}{\lambda} z_{n}  \\
             &=& \frac{8}{\lambda}\sum_{i=2}^{n-1}z_{i}.
\end{eqnarray*}
\noindent
The equality follows, since $v_1$ and $v_n\in\partial\mathbb{R}^3_+=\mathbb{R}^2$. $\square$\\

\begin{ex}
The cubical trefoil $1$-knot of this example was studied by Yunnan Diao in \cite{diao} (see Figure \ref{trefoildiao}). It consists of  24 edges and  13 corner vertices, so the cubical arc  $A_{1}$ has 
22 vertices $v_{1}$, $v_{2}$, $\ldots$, $v_{22}$.  Hence the cubical spun trefoil $2$-knot $C\Spin(A_1)$ obtained from  $A_{1}$ has area
$\Area(C\Spin(A_1))=288$. For this example, we have that $\Long(A_{1})=21$.\\

\begin{figure}[h] 
\begin{center}
\includegraphics[width=6cm, height=6cm ]{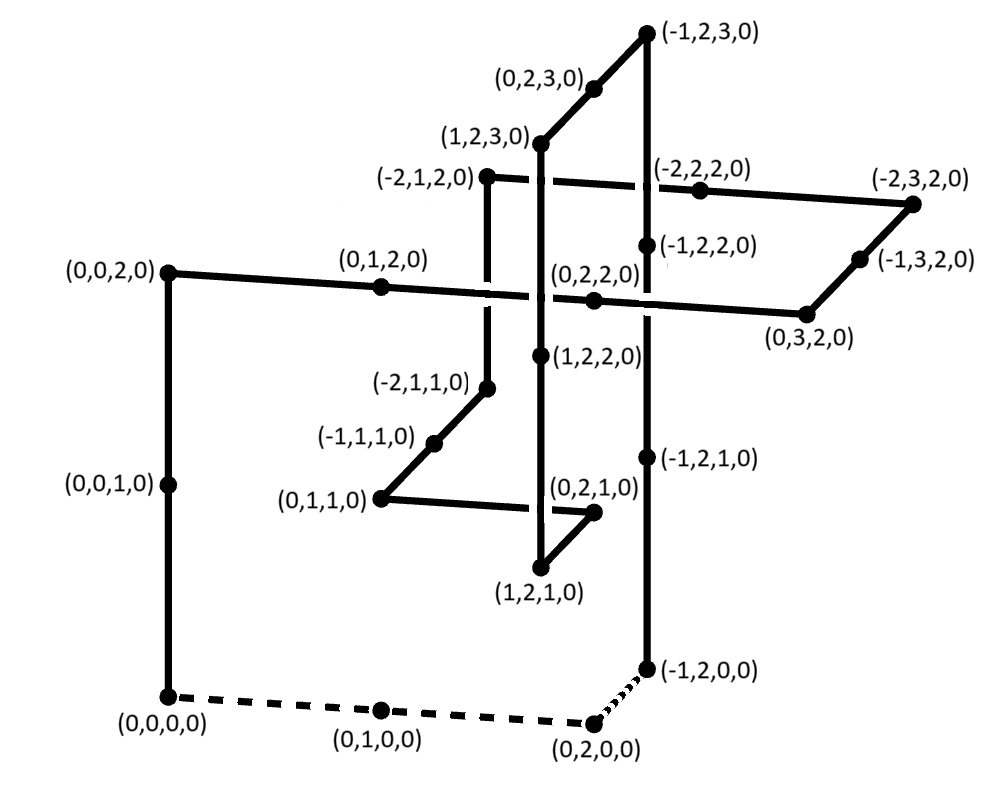}
\caption{\sl Cubical trefoil knot such that $\Area(C\Spin(A_1))$ is 288.}
\label{trefoildiao}
\end{center}
\end{figure}
\end{ex}

\begin{ex}
The cubical trefoil $1$-knot of this example was studied by A. Baray and G. Hinojosa in \cite{BH} (see Figure \ref{trefoilgaby}). It has 24 edges and 14 corner vertices. Then, the cubical arc $A_{2}$ is given by the 22 vertices $v_{1}$, $v_{2}$, $\ldots$, $v_{22}$. The cubical spun trefoil $2$-knot $C\Spin(A_2)$ obtained from  
$A_{2}$ has area $\Area(C\Spin(A_2))=256$. For this example, we have that $\Long(A_{2})=21$.

\begin{figure}[h] 
\begin{center}
\includegraphics[width=6cm, height=6cm]{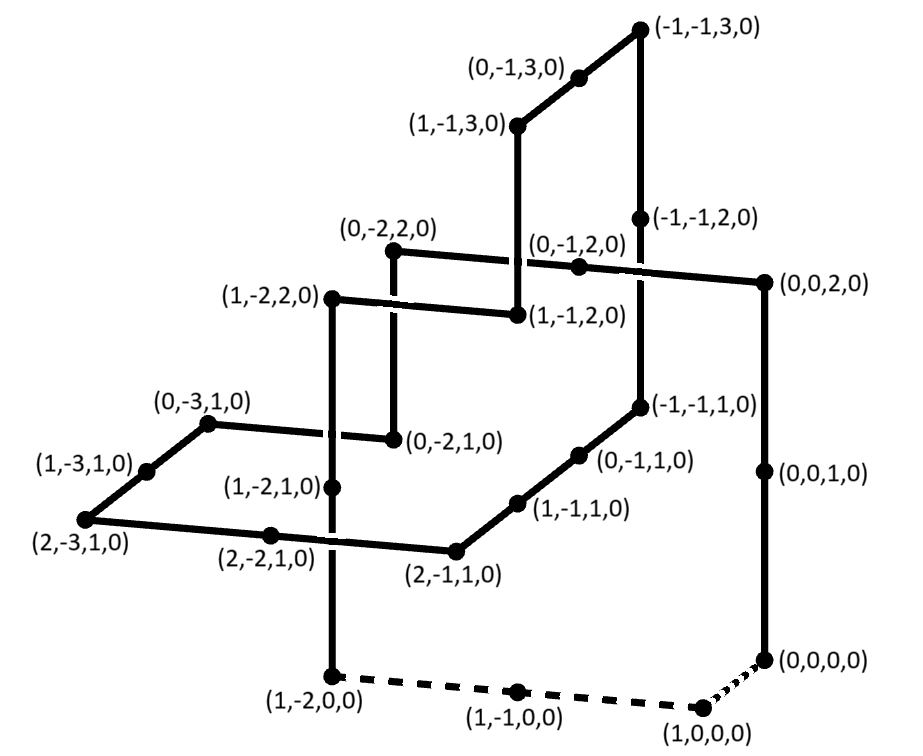}
\caption{\sl Cubical trefoil knot such that $\Area(C\Spin(A_2))$ is 256.}
\label{trefoilgaby}
\end{center}
\end{figure}
\end{ex}
\vskip .3cm

\begin{ex}\label{ExJGR}
The cubical trefoil $1$-knot of this example was given by the authors. It consists of  26 edges and  13 corner vertices, and the cubical arc $\mathcal{A}$ has 
23 vertices $v_{1}$, $v_{2}$, $\dots$, $v_{23}$.  Hence the cubical spun trefoil $2$-knot $C\Spin({\cal{A}})$ obtained from $\cal{A}$ has area
$\Area(C\Spin({\cal{A}}))=216$. For this example, we have that $\Long({\cal{A}})=22$, see Figure \ref{trefoiljj}.\\

\begin{figure}[h] 
\begin{center}
\includegraphics[width=7cm, height=7cm ]{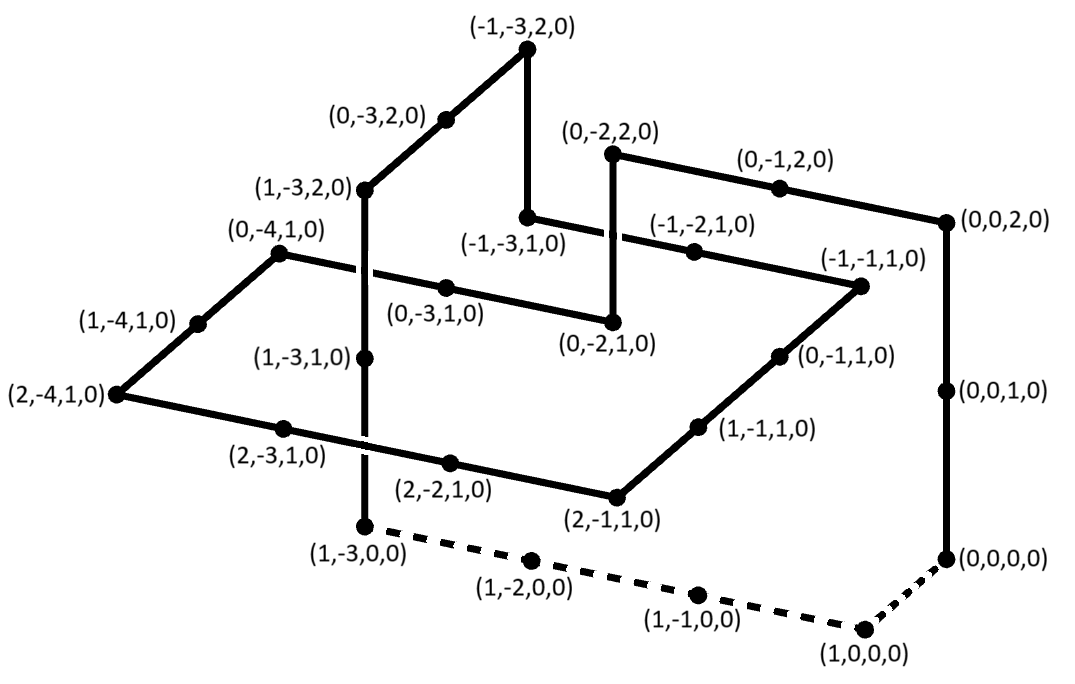}
\caption{\sl Cubical trefoil knot such that $\Area(C\Spin({\cal{A}}))$ is  216.}
\label{trefoiljj}
\end{center}
\end{figure}

\end{ex}

\subsection{Area of reduced cubical spun 2-knots}

We are interested in the minimal area of cubical spun 2-knots, so we analyze whenever a subcubical arc $A$ can yield a cubical spun 2-knot.\\

\noindent Let $A$ be a subcubical arc contained in ${\cal{C}}_{\lambda}$ with vertices  $v_{1}$, $v_{2}$, $\ldots$, $v_{n}$ in $\mathbb{R}_{+}^{3}$, and edges $\overline{v_{1}v_{2}}$, $\overline{v_{2}v_{3}}$, $\ldots$, $\overline{v_{n-1}v_{n}}$, such that $v_{1}$ and  $v_{n}$ are in  $\mathbb{R}^{2}$, and $v_{2}$, $v_{3}$, $\ldots$, $v_{n-1}$ are in the interior of $\mathbb{R}_{+}^{3}$. Hence $\Spin\big(A\big)=\cup_{i=1}^{n-1}\Spin(\overline{v_{i}v_{i+1}})$ is a $2$-sphere embedded in $\mathbb{R}^{4}$ obtained by gluing disks, annulus, and cylindrical surfaces in a proper way, where $v_{i+1}-v_{i}=\frac{1}{\lambda} e_{k_{i}}$, for $k_{i} \in \{ \pm1, \pm2, \pm3 \}$, $\frac{1}{\lambda}=||v_{i+1}-v_{i}||$. Notice that the surface $\Spin(\overline{v_{i}v_{i+1}})$ can be deformed to get a subcubical surface, but in general, the coordinates of its vertices are not integers. We would like to translate the subcubical spun knot $C\Spin\big(A\big)_{\lambda}$ by a vector $b$ in such a way that $b+C\Spin\big(A\big)_{\lambda}$ is
contained in the 2-skeleton of the canonical cubulation ${\cal{C}}$ of $\mathbb{R}^4$. This is possible only for $\lambda=2$, since the spin of each vertex
$v_{i}=(x_{i},y_{i},z_{i},0)\in A$ is given by  $\Spin(v_i)=(x_i, y_i, z_{i}\cos (2\pi \theta), z_i\sin(2\pi \theta))$, we have that $x_i$ and $y_i$ must be integers, and 
$z_i$ must be a multiple of $\frac{1}{2}$, because the diameter of the corresponding surface $\Spin(\overline{v_{i}v_{i+1}})$ which is equal to $2z_i$ is an integer. Thus, we require that  if 
$\alpha_{i}=||v_{i+1}-v_{i}||$, then

\begin{enumerate}
\item $\alpha_{1}=\frac{1}{2}$.
	\item $\alpha_{i}=1$ para $i=2,3,\ldots, n-2$.
	\item $\alpha_{n-1}=\frac{1}{2}$.
\end{enumerate}
\noindent In this case, the translation vector is $b=(0,0,\frac{1}{2},\frac{1}{2})$.\\

\noindent A subcubical arc $A$ with the $\alpha_i$'s satisfying the previous conditions 
will be called a {\it reduce cubical arc}. Observe that if $A$ is a reduce cubical arc, it is still a subcubical arc contained in ${\cal{C}}_{2}$; however,  our description omits some vertices of ${\cal{C}}_{2}$, since we are interested in the cubical 2-knot  $b+C\Spin\big(A\big)_{2}$, hence we consider only those vertices that are sent under this translation to the lattice $\mathbb{Z}^4$ of $\mathbb{R}^4$.\\

\begin{definition}
We say that a  $2$-knot is a reduced cubical spun $2$-knot if it is obtained from a reduced cubical arc $A$ via the above process. It will be denoted by $\RCSpin(A)$, that is, $RC\Spin(A)=b+C\Spin(A)_{2}$.\\
\end{definition}

\noindent Using the previous notation, we can calculate the area of the reduced cubical spun knot $RC\Spin(A)$ by Theorem 2.

\begin{coro}\label{areaR}
Let $RC\Spin(A)$ be a reduced cubical spun $2$-knot  obtained from the reduced cubical arc  $A$. Let $v_{1}$, $v_{2}$, $\ldots$, $v_{n}$ be the vertices of $A$, where $v_{i}=(x_{i},y_{i},z_{i},0)$, such that $\alpha_{i}=||v_{i+1}-v_{i}||$ satisfies that $\alpha_{1}=\frac{1}{2}$, $\alpha_{i}=1$ for $i=2,3,\ldots, n-2$ and  $\alpha_{n-1}=\frac{1}{2}$. Then 
$$\Area(RC\Spin(A))=\left(8\sum_{i=2}^{n-1}z_{i}\right) -2. $$
\end{coro}

\noindent {\it Proof.} By an argument similar to that used in Theorem 2, we have that  
\begin{eqnarray*} 
\Area(RC\Spin(A)) &=& 4\sum_{i=1}^{n-1}\alpha_{i}(z_{i}+z_{i+1}) \\
             &=& 4\alpha_{1}(z_{1}+z_{2}) + 4\sum_{i=2}^{n-2}\alpha_{i}(z_{i}+z_{i+1})+4\alpha_{n-1}(z_{n-1}+z_{n}).
\end{eqnarray*}
Since $\alpha_{1}=\frac{1}{2}$, $\alpha_{i}=1$ for $i=2,3,\ldots, n-2$ and $\alpha_{n-1}=\frac{1}{2}$, then 
\begin{eqnarray*} 
\Area(RC\Spin(A)) &=& 4\left(\frac{1}{2}\right)(z_{1}+z_{2}) + 4\sum_{i=2}^{n-2}(z_{i}+z_{i+1})+4\left(\frac{1}{2}\right)(z_{n-1}+z_{n}) \\
             &=& 2(z_{1}+z_{2}) + 4\sum_{i=2}^{n-2}(z_{i}+z_{i+1})+2(z_{n-1}+z_{n}) \\
             &=&  4\sum_{i=1}^{n-1}(z_{i}+z_{i+1}) -2(z_{1}+z_{2})-2(z_{n-1}+z_{n}). 
\end{eqnarray*}
\noindent By  Theorem 2, we have that 
 $4\sum_{i=1}^{n-1}(z_{i}+z_{i+1})= 8\sum_{i=2}^{n-1}z_{i}$, 
but $z_{1}=z_{n}=0$ y $z_{2}=z_{n-1}=\frac{1}{2}$, hence
$$
\Area(RC\Spin(A))=\left(8\sum_{i=2}^{n-1}z_{i}\right) -2. 
$$\hfill $\square$\\

\begin{rem} \label{semicub}
Given a cubical arc $A$ with vertices $v_{1}$, $v_{2}$, $\ldots$, $v_{n}$, where each $v_i$ belongs to the lattice $\mathbb{Z}^4$, each edge has length one and is parallel to some coordinate axis, we can obtain a reduced cubical arc $A'$ by modifying the $z$-coordinates of the $v_i$'s that are different from zero, such that the vertices of $A'$ are $v'_{1}=v_1$, $v'_{2}$, $\ldots$, $v'_{n-1}$, $v'_n = v_n \in(\frac{1}{2}\mathbb{Z})^4$, where  
$v'_{i}= (x_{i},y_{i},z'_{i},0)=(x_{i},y_{i},z_{i}-\frac{1}{2},0)$, for $i=2, 3, \dots, n-1$. 
Thus, $A'$ is a reduced cubical arc so that we can establish a one-to-one correspondence between the set of cubical arcs and the set of 
reduced cubical arcs with the properties just described.
\end{rem}

\vspace{.1in}

\begin{coro} \label{areared}
Let $C\Spin(A)$ be a  \emph{cubical spun $2$-knot} obtained from the cubical arc  $A$ with vertices $v_{1}$, $v_{2}$, $\dots$, $v_{n}\in\mathbb{Z}^4$, where $v_{i}=(x_{i},y_{i},z_{i},0)$. Then, the minimal area of $C\Spin(A)$ is upper bounded by $\displaystyle{8 \left(\sum_{i=2}^{n-1}z_{i} \right) -4n+6}$; that is,
$$
A(C\Spin(A))\leq 8 \left( \sum_{i=2}^{n-1}z_{i} \right) -4n+6.
$$
\end{coro}

\noindent {\it Proof.} 
Let $A'$ be the cubical arc obtained from $A$, by modifying the $z$-coordinate of each of its vertices $v_2, v_3, \dots, v_{n-1}$, such that the vertices of $A'$ are
$v_{1}$, $v'_{2}$, $\ldots$, $v'_{n-1}$, $v_n \in(\frac{1}{2}\mathbb{Z})^4$, where  $v'_{i}= (x_{i},y_{i},z'_{i},0)=(x_{i},y_{i},z_{i}-\frac{1}{2},0)$, for $i=2,3, \dots, n-1$. Then, $A'$ is a reduced cubical arc, so that we can consider  
the corresponding
reduced cubical spun $2$-knot $RC\Spin(A')$. By Corollary 1, 
\begin{eqnarray*} 
\Area(RC\Spin(A')) &=& \left(8\sum_{i=2}^{n-1}z'_{i}\right) -2\\
             &=& \left(8\sum_{i=2}^{n-1}\left(z_{i}-\frac{1}{2}\right)\right) -2\\
             &=& \left(8\sum_{i=2}^{n-1}z_{i}\right)-4n +6.\\
\end{eqnarray*}
\noindent Since $RC\Spin(A')$ is isotopic to $C\Spin(A)$, the result follows. $\square$\\

\begin{ex}\label{RExJGR}
Consider the cubical $1$-knot described in example \ref{ExJGR}. So, it consists of  26 edges and  13 corner vertices, and the cubical arc $\mathcal{A}$ has 
23 vertices $v_{1}$, $v_{2}$, $\dots$, $v_{23}$. Now, we modify the cubical arc $\cal{A}$ to get a semi-subcubical arc ${\cal{A}}'$ Hence the reduced cubical spin trefoil $2$-knot $RC\Spin({\cal{A}}')$ obtained from $\cal{A}'$ has 
$\Area(RC\Spin({\cal{A}}'))=130$, see Figure \ref{trefoilreduced}.
\begin{figure}[h] 
\begin{center}
\includegraphics[width=7cm, height=7cm ]{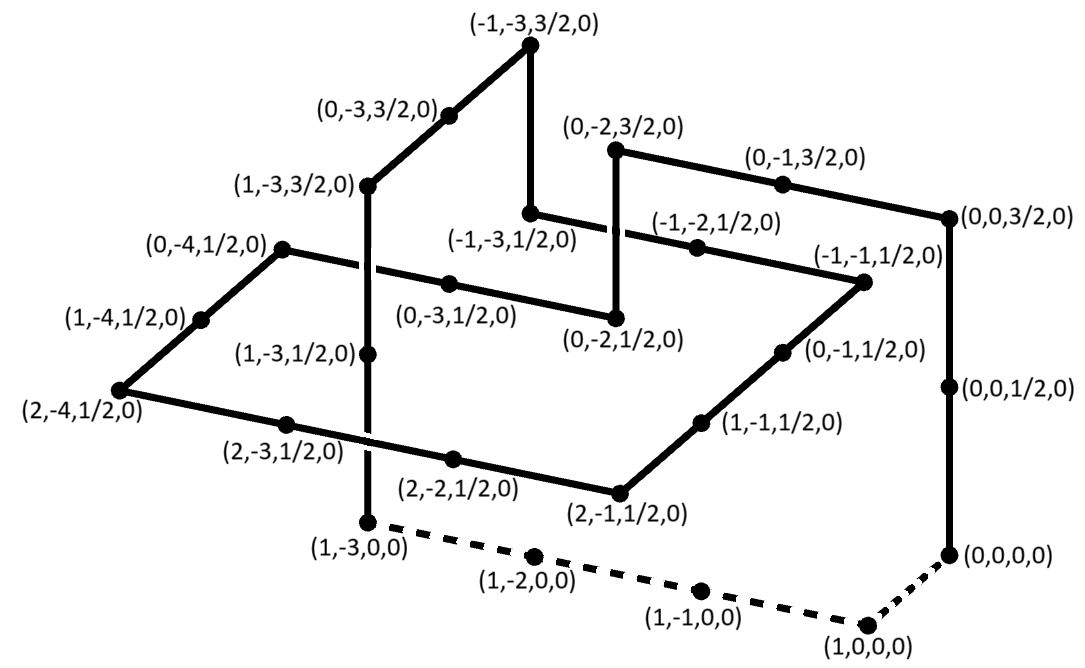}
\caption{\sl Reduced cubical trefoil arc such that $\Area(RC\Spin({\cal{A}}'))$ is 130.}
\label{trefoilreduced}
\end{center}
\end{figure}

\end{ex}

\subsection{Cubical spin minimal area}

\noindent  We denote the set of all cubical spun $2$-knots and reduced cubical spun $2$-knots  by ${\mathcal{CS}}^2$; \emph{i.e.},
$$
{\mathcal{CS}}^2=\{C\Spin(A)\,:\,\mbox{$A$ is a cubical arc}\}\cup \{RC\Spin(A')\,:\,\mbox{$A'$ is a reduced cubical arc}\}.
$$

\noindent Consider the subset ${\cal{CS}}(K^2)\subset {\mathcal{CS}}^2$ consisting of all cubical spun 2-knots and reduced cubical spun 2-knots isotopic to 
the spun 2-knot $K^2$, that is
$$
{\cal{CS}}(K^2)=\{\tilde{K}^2\in {\mathcal{CS}}^2\,\,|\,\,\tilde{K}^2\,\,\mbox{  is isotopic to}\,\,K^2\}.
$$

\begin{definition}
The \emph{cubical spin minimal area} of a spun 2-knot $K^{2}$, denoted by  $A_{\mathcal{CS}}(K^{2})$, is defined as the minimal area of all cubical spun 2-knots isotopic to $K^{2}$, \emph{i.e.,}
$$
A_{\mathcal{CS}}(K^{2})=\min\{\Area(\tilde{K}^{2})\,:\,\tilde{K}^2\in {\cal{CS}}(K^2)\}.\\
$$
\end{definition}

\begin{rem}
 Notice that the cubical spin minimal area of spun knots is also an invariant in the set ${\mathcal{CS}}^2$. In particular, the cubical spin minimal area of a spun 2-knot $K^2$ is greater than or equal to the minimal area of $K^2$; \emph{i.e.}, 
$A(K^{2})\leq A_{\mathcal{CS}}(K^{2})$.\\
\end{rem}

\section{Minimal area of cubical spin trefoil 2-knot}

In this section, we will start computing the cubical spin minimal area of the spun trefoil 2-knot that is realised by the reduced cubical spin $2$-knot $RC\Spin(A')$ given in example 5. \\


\begin{main3}\label{main3}
The cubical spin minimal area of the spun trefoil knot $\Spin(T_{2,3})$ is 130, that is, 
$$
A_{\mathcal{CS}}(\Spin(T_{2,3}))=130.
$$
\end{main3}
\vskip .25cm
\noindent Before proceeding to prove Theorem 3, let us give some general considerations about any cubical 1-knot $K$ isotopic to the trefoil knot $T_{2,3}$. Suppose that $K$ has $m$ vertices, so let the vertices of $K$ be $v_i=(x_i, y_i, z_i,0)$, for $i=1, 2, \dots, m$.\\

\begin{rem} \label{rem1}
First of all, if all the vertices of $K$ have coordinates in $z$ equal to 0 or 1, it follows that $K$ is unknotted; hence, there are at least two vertices of $K$ with coordinates in $z$ equal to 2. Now suppose $K$ cannot be reduced in length by a cubulated move, then if there are only two vertices in $K$ with coordinates in $z$ equal to 2, they are consecutive vertices, say $v_i$ and $v_{i+1}$, and they should form an edge $e_i$ of $K$, so there cannot be an edge of $K$ directly below $e_i$, since each of these vertices is joint down to a vertex with coordinate in $z$ equal to 1 by edges $e_{i-1}$ and $e_{i+1}$, respectively, then we can ``replace" the arc $e_{i-1}e_ie_{i+1}$ applying a cubulated move of (M2)-type 
and still obtain a knot isotopic to $T_{2,3}$, but with a smaller length.\\
\end{rem}

\begin{figure}[h] 
\begin{center}
\includegraphics[width=8cm, height=4cm ]{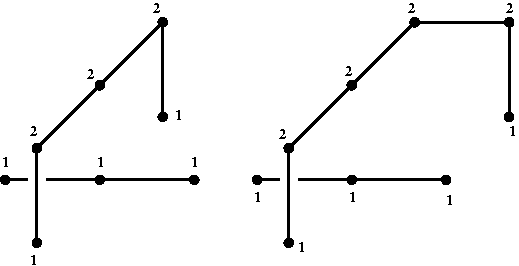}
\caption{\sl Two examples of crossing steps.}
\label{bridges}
\end{center}
\end{figure}

\noindent By the previous Remark \ref{rem1}, if one vertex of $K$ has coordinate in $z$ equal to 2, there are three (or more) consecutive vertices $v_i$, $v_{i+1}$, $v_{i+2}$ of $K$ with coordinate in $z$ equal to 2. We define a ``crossing step" of $K$ as a series of consecutive vertices $v_{i}, v_{i+1}, \dots, v_{i+k-1}$, with coordinate in $z$ equal to 2, for some $k \geq 3$ together with the vertices $v_{i-1}$ and $v_{i+k}$ with coordinate in $z$ equal to one, and three consecutive vertices with coordinate in $z$ equal to 1, directly below the crossing step, making a crossing of $K$, see the Figure \ref{bridges} for examples of two crossing steps, where the numbers indicate the $z$ coordinate in the corresponding vertex. Observe that a crossing step consists of at least eight vertices in total.\\

\begin{figure}[h] 
\begin{center}
\includegraphics[width=6cm, height=5cm ]{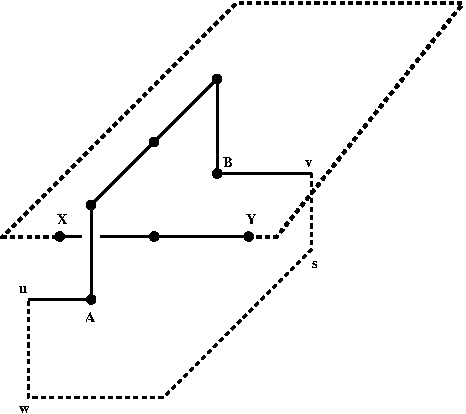}
\caption{\sl Two disjoint arcs from a crossing step.}
\label{ejemplo}
\end{center}
\end{figure}

\begin{rem} \label{rem2}
If $K$ has only one crossing step and the other vertices have coordinates in z equal to 1 or 0, then $K$ is unknotted.  To see this, observe that a crossing step consists of two disjoint arcs $AB$ and $XY$ of $K$ (see Figure \ref{ejemplo}). We know that there are two pair of vertices of $K$, say $u=(x_i, y_i, 1, 0)$, $v=(x_j, y_j, 1, 0)$ and $w=(x_i, y_i,0,0)$, $s=(x_j, y_j,0,0)$ such that $u$ is connected down directly to $w$
and $v$ is connected down directly to $s$. The four ends of the two disjoint arcs $AB$ and $XY$, that is, the vertices $A$, $B$, $X$ and $Y$ have to be connected between them or to $u$ or $v$ (maybe travelling before in some edges with coordinate in $z$ equal to 1).
If $A$ and $B$ are connected to $u$ and $v$ in some order, then $X$ and $Y$ are connected themselves (as shown in Figure \ref{ejemplo}); this implies that $K$ is a link, which is not possible. In the same way, if $X$ and $Y$ are connected to $u$ and $v$, then $A$ and $B$ are connected between them, and then $K$ consists of two disjoint closed arcs, which is a contradiction. Therefore, one of the ends $A$ or $B$ should be connected to one of the ends $u$ or $v$, say $u$, then one of the ends $X$ or $Y$ must be connected to $v$. In each of these last cases, we obtain that $K$ is unknotted, since it will have only one cross.\\
\end{rem}

\noindent In other words, for $K$ to be a knot with vertices with coordinates in $z$ equal to 0, 1, and 2, it must be at least two crossing steps, that is, at least two crosses. We are ready to prove Theorem 3.\\

\begin{main3}\label{main3}
The cubical spin minimal area of the spun trefoil knot $\Spin(T_{2,3})$ is 130, that is, 
$$
A_{\mathcal{CS}}(\Spin(T_{2,3}))=130.
$$
\end{main3}
\vskip .25cm
\noindent {\it Proof.} Let $K$ be a cubical 1-knot isotopic to the trefoil knot $T_{2,3}$, and let $K^2$ the corresponding cubical spun 2-knot obtained from the cubical arc $A$ generated by $K$, that is, $K^2$ is a
cubical spun 2-knot isotopic to the spun trefoil knot $\Spin(T_{2,3})$. 
Let $A'$ be the cubical arc obtained from $A$ by modifying the $z$-coordinate of each of its vertices, that is, the vertices of $A'$ are
$v'_{1}=v_1$, $v'_{2}$, $\dots$, $v'_{n-1}$, $v'_n=v_n \in(\frac{1}{2}\mathbb{Z})^4$, where  $v'_{i}=(x_{i},y_{i},z_{i}-\frac{1}{2},0)$, 
for $i=2, 3, \dots, n-1$,
then  
$A'$ is a reduced cubical arc.\\

\noindent Now, consider the corresponding
reduced cubical spin $2$-knot $RC\Spin(A')$, by Corollary \ref{areared}, 
$\Area(RC\Spin(A'))= \displaystyle{8 \left(\sum_{i=2}^{n-1}z_{i} \right) -4n+6}$.
Let us assume that 
$\Area(RC\Spin(A')) = \displaystyle{ 8 \left(\sum_{i=2}^{n-1}z_{i} \right) -4n+6  < 130}$,
then, by Theorem 2 making $\lambda=1$, we have that $\displaystyle{\Area(K^2)= 8 \sum_{i=2}^{n-1}z_{i} < 124 + 4n}$. We will prove that $\Area(K^2)$ can not be less than 216. 
This will imply that $216 \leq \Area(K^2) < 124 +4n$, then $23 < n$ or
$24 \leq n$. Since $K$ is a knot it must have at least six coordinates in $z$ equal to 2 and it has in total $n-2$ coordinates in $z$ distinct from zero, so $\displaystyle{8 \left(\sum_{i=2}^{n-1}z_{i} \right)} -4n+6 \geq 8(6(2) + ((n-2)-6)) -4n +6= 8(n+4) -4n + 6 = 4n+38 \geq 4(24)+38=134$, which is a contradiction.
Now we can take $n=23$ as in Example \ref{RExJGR}, reaching the minimal area $\Area(RC\Spin(Arc({\cal{A}})))=130$, and since
the cubical knots in Examples \ref{ExJGR} and \ref{RExJGR} are isotopic we will be done.

\vspace{.3in}

\noindent Let us prove now that $\Area(K^2) < 216$ is impossible. Assume that $K$ has $V_K=m$ vertices $v_i=(x_i, y_i, z_i, 0)$, with $z_i \neq 0$, for $i=1, 2, \dots, m$, so the total number of vertices of $K$ is $V_K$ plus the number of vertices with coordinate in $z$ equal to zero, then
$\displaystyle{\Area(K^2)=8\sum_{i=1}^{m}z_{i} < 216}$, implies that $0 < \displaystyle{\sum_{i=1}^{m}z_{i} < 27}$,
that is $\displaystyle{\sum_{i=1}^{m}z_{i} \leq 26}$. Set $S(m)=\displaystyle{\sum_{i=1}^{m}z_{i}}$, we will proceed by cases depending on the number of $z_i$'s different from zero to get a contradiction. \\

\noindent Assume first that $V_K=m \geq 26$, then $26 \geq S_m \geq 26$, thus the $m$ vertices satisfies that $z_i =1$, which is impossible, otherwise $k$ is unknotted. Then $V_K=m \leq 25$.\\

\noindent Case $V_K= 25$. Suppose there is some vertex with $z_i=3$ (or greater than 3), then there are at least two vertices with coordinate $z$ equal to 2, and it follows that 
$S_{25} \geq 3(1)+2(2)+22 > 26$, a contradiction.
Therefore, $1 \leq z_i \leq 2$ for all $i$, and the same will hold, by the same estimation, for $V_K \leq 24$. If there is a vertex with coordinate $z_i=2$, there are at least another two vertices with the coordinate in $z$ equal to $2$ by Remark \ref{rem1}, then 
$S_{25} \geq 3(2)+22 > 26$, once again a contradiction, 
therefore in this case all $z_i=1$, which is not possible.\\

\noindent Case $V_K= 24, 23$. 
When there is a vertex with $z_i=2$, then
$S_m \geq 3(2)+(m-3)=m+3 > 26$ for $m=24$ and we arrive to the same contradiction. When $m=23$, we get that $S_{23}$ is exactly 26, so $K$ consists of three vertices with coordinate in $z$ equal to two, and the rest have coordinate in $z$ equal to 1 (or equal to 0), by Remark \ref{rem2} this is not possible.\\

\noindent Case $V_K= 22, 21$. If there are three vertices with coordinates in $z$ equal to 2, then 
$S_{m} \geq 3(2)+(m-3)=m+3$, for which $S_{22}$ is equal to 25 or 26, and
$S_{21}$ is equal to 24, 25, or 26. Then, either  
$K$ has three, four or five consecutive vertices with coordinate in $z$ equal to 2 (remember that we cannot have an isolated vertex with coordinate in $z$ equal to 2 or only two consecutive vertices with coordinate in $z$ equal to 2 by Remark \ref{rem1}), and the other vertices with coordinates in $z$ equal to 1 and 0, by Remark \ref{rem2} this is a contradiction. Therefore, the other possibility is to have
two different crossing steps, so  $S_{21} \geq 6(2)+(15)=27 > 26$, which is a contradiction.\\

\noindent Case $V_K= 20$.  In this case $S_{20} \geq 6(2)+14=26$, since we must have two different crossing steps, and observe that $S_{20}$ is exactly equal to 26, that is, $K$ has six vertices with coordinate in $z$ equal to 2 forming two different crossing steps, and 14 vertices with coordinate in $z$ equal to 1.
Since each of these crossing steps should consist of 8 vertices, let us see how these crossing steps can be; in fact, there are three possibilities.\\

\begin{figure}[h] 
\begin{center}
\includegraphics[width=8cm, height=6cm ]{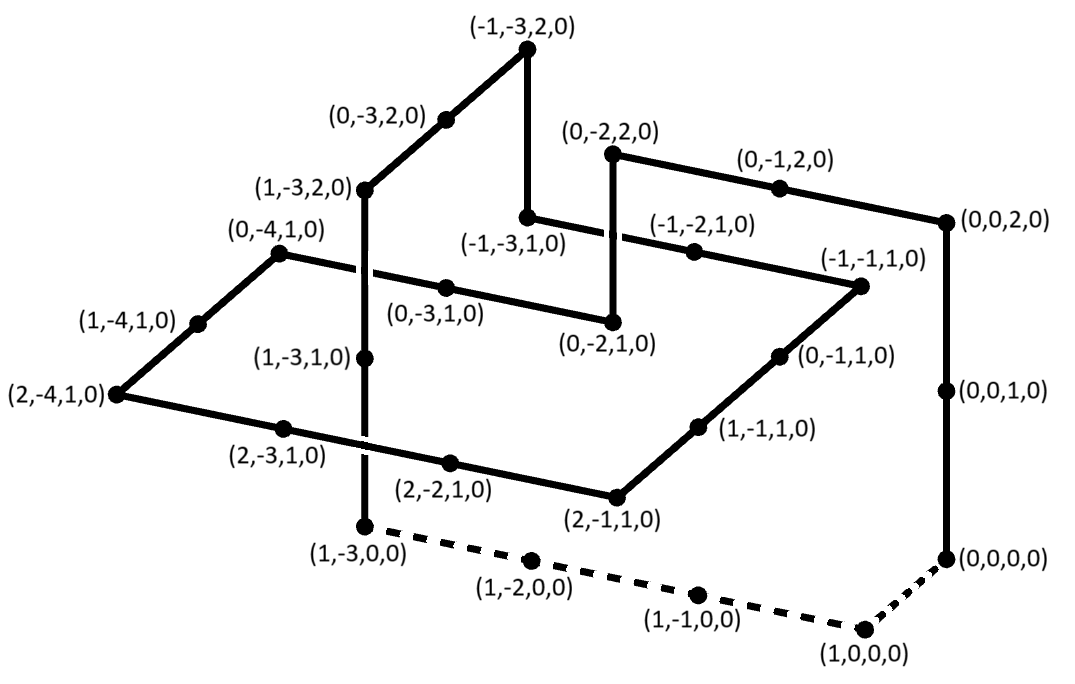}
\caption{\sl Two crossing steps using only 14 vertices.}
\label{trefoilr}
\end{center}
\end{figure}

\begin{figure}[h] 
\begin{center}
\includegraphics[width=8cm, height=6cm ]{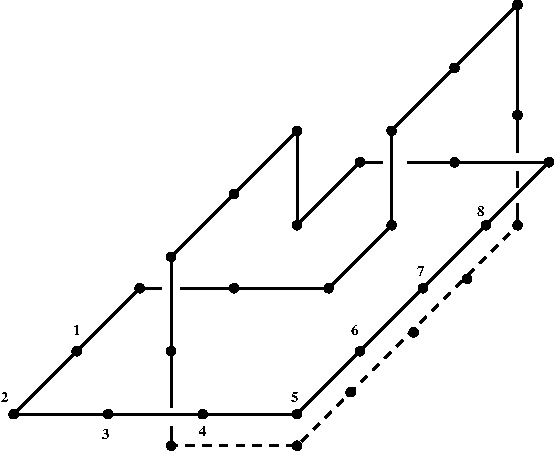}
\caption{\sl Two crossing steps using 16 vertices.}
\label{trefoilm}
\end{center}
\end{figure}

\noindent In Figure \ref{trefoiljj}, we see an example of a knot with two crossing steps sharing a vertex between them and using 21 vertices with coordinates in $z$ different from zero. Meanwhile in Figure \ref{trefoilr}, we present an example with two crossing steps sharing two vertices, so between the two crossing steps there are 14 vertices (this is the minimal number of vertices needed to construct two different crossing steps), and $K$ has in total 21 vertices with coordinate in $z$ different from zero. Any other configuration of $K$ with two crossing steps will use at least 16 vertices, see Figure \ref{trefoilm}, where the configuration have two disjoint crossing steps (from Figure \ref{trefoilr}, move one of the crossing steps 1 unit in the $xy$-plane), but in total there are 24 vertices with coordinate in $z$ different from zero. Observe that in the three previous examples, we have used the minimal number of edges to close the knot once the crossing steps are constructed.\\

\noindent In the three examples of the knot with two crossing steps, we observe that the crossing steps generate certain total width for $K$ in the $x$ direction and the $y$ direction: in Figure \ref{trefoiljj} we have a total width of 2 units in the $x$ direction and four units in the $y$ direction, in Figure \ref{trefoilr} we have total width of 3 units in both $x$ and $y$ directions, and in Figure \ref{trefoilm} we have a total width of 4 units in the $x$ direction and three units in the $y$ direction, and they use 15, 14 and 16 total vertices to form the two crossing steps, respectively. \\

\noindent By the same argument in Remark \ref{rem2}, one end of each crossing step should be connected (maybe after traveling by some edges in the $xy$-plane at height 1) to a vertex of $K$ with coordinate in $z$ equal to zero. In this way, the other ends of each crossing step should be connected. In each possible configuration of the crossing steps, we will determine how many vertices with coordinates in $z$ equal to one are required. Remember that $S_{20}=26$ and $K$ has 6 vertices with coordinate in $z$ equal to 2, and 14 vertices with coordinate in $z$ equal to 1.\\

\noindent In the case when the two crossing steps share two vertices, as in Figure \ref{trefoilr}, the two crossing steps already use 14 vertices, so there are six vertices left with coordinate in $z$ equal to 1. Since we have total width of 3 units in both $x$ and $y$ directions, the minimal number of edges in the $xy$-plane at height $z=1$ to go around the end of the crossing step that is connected to the $z=0$ plane is $8$, so we need seven vertices more to closed the knot (in the Figure \ref{trefoilr} the vertices are $(-1,0,1,0)$, $(-1,1,1,0)$, $(0,1,1,0)$, $(1,1,1,0)$, $(2,1,1,0)$, $(2,0,1,0)$, $(2,-1,1,0)$), which is a contradiction since we can use only 20 vertices. \\

\noindent When the two crossing steps share only one vertex, as in Figure \ref{trefoiljj}, we have used 15 vertices to form the two crossing steps, so we have five vertices left with coordinate in $z$ equal to 1. But since we have a total width of 2 units in the $x$ direction and four units in the $y$ direction, there are seven edges and six vertices needed to closed the knot (in the Figure \ref{trefoiljj} the vertices are $(-1,-2,1,0)$, $(2,-1,1,0)$, $(2,-2,1,0)$, $(2,-3,1,0)$, $(2,-4,1,0)$, $(1,-4,1,0)$), and we get a total of 21 vertices which is a contradiction.\\

\noindent Finally, for the case when the crossing steps use 16 vertices (so there are only four vertices that are not in the crossing steps) as in Figure \ref{trefoilm}, we see that we need at least eight vertices (the ones marked with the numbers in Figure \ref{trefoilm}) to complete the knot, which is more than 20 vertices in total. This completes the proof, since it is impossible to have a knot with fewer than 20 vertices and coordinates in $z$ equal to 1, containing two different crossing steps, as previously argued. $\square$\\

\subsection{Main Theorem}

In this section, we will prove our main theorem. We will start describing with more detail our reduced cubical spun knot $RC\Spin({\cal{A'}})$ (see Figure \ref{trefoilreduced} in Example 5). Let  $v_{1}$, $v_{2}$, $\ldots$, $v_{23}$ be the vertices of  ${\cal{A'}}$ and let $\overline{v_{1}v_{2}}$, $\overline{v_{2}v_{3}}$, $\ldots$, $\overline{v_{22}v_{23}}$ be its edges.  In the following, we will describe $RC\Spin({\cal{A'}})$ as the union of eleven cubical surfaces, where each of them is obtained spinning one or more consecutive edges of ${\cal{A'}}$ according to the construction given in Section 3.1.  \\

\begin{figure}
\centering
\includegraphics[height=11.5cm]{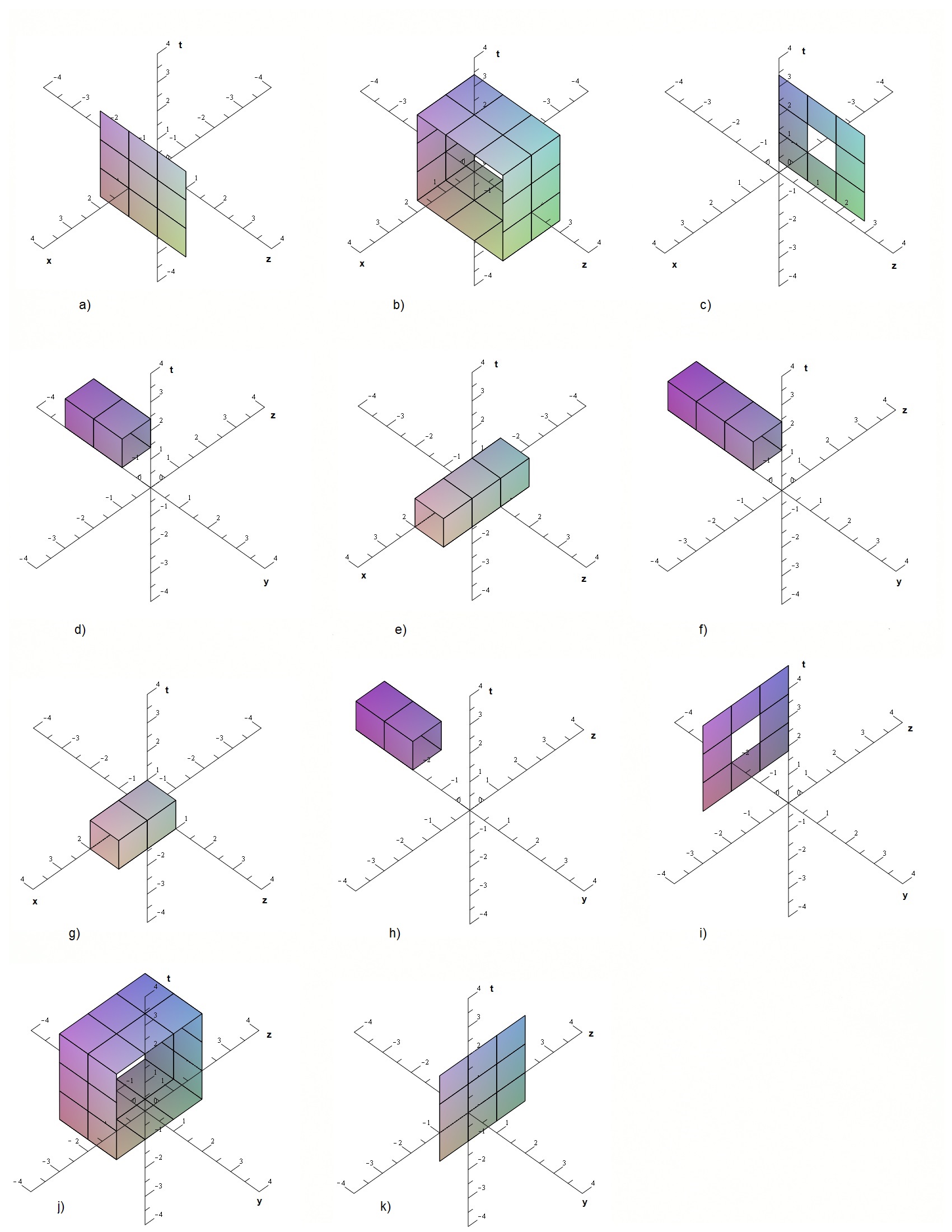}
\caption{The reduced cubical spun knot $RC\Spin({\cal{A'}})$ decomposed into eleven cubical surfaces.}
\label{piece}
\end{figure}

\noindent Step 1. The square $D_{1}$ is obtained from the union of the two edges $\overline{v_{1}v_{2}}\cup\overline{v_{2}v_{3}}$. Its length side is three (see Figure \ref{piece} part $a$).  So, its area is 9.\\

\noindent Step 2. The square cylinder $D_{2}$ is gotten from $\overline{v_{3}v_{4}}\cup\overline{v_{4}v_{5}}$. Its base length is three and its height is two (see Figure \ref{piece} part $b$). Hence, its area is 24.\\

\noindent Step 3. The square annulus $D_{3}$ comes from the edge $\overline{v_{5}v_{6}}$. Its largest side length is three, and the shortest side length is 1 (see Figure \ref{piece} part $c$). Thus, its area is 8. \\

\noindent Step 4. The square cylinder $D_{4}$ is obtained from $\overline{v_{6}v_{7}}\cup\overline{v_{7}v_{8}}$ Its base length is one with height two (Figure \ref{piece}  part $d$). Therefore, its area is 8.\\

\noindent Step 5. The square cylinder $D_{5}$ is constructed via the edges $\overline{v_{8}v_{9}}$, $\overline{v_{9}v_{10}}$, and $\overline{v_{10}v_{11}}$. Its base length is one
and its height is three (Figure \ref{piece} part $e$). So, its area is 12.\\

\noindent Step 6. The square cylinder $D_{6}$ is comes from  $\overline{v_{11}v_{12}}\cup \overline{v_{12}v_{13}}\cup \overline{v_{13}v_{14}}$. Its base length is one with height three (Figure \ref{piece}  part $f$). Hence, its area is 12.\\

\noindent Step 7. The square cylinder $D_{7}$ is gotten from  $\overline{v_{14}v_{15}}$ and  $\overline{v_{15}v_{16}}$. Its base length is one, with a height of two (see Figure \ref{piece} part $g$). This implies that its area is 8.\\

\noindent Step 8. The square cylinder $D_{8}$ is constructed via the edges $\overline{v_{16}v_{17}}$ and $\overline{v_{17}v_{18}}$. Its base length is one and its height is two (Figure \ref{piece} part $h$). Thus, its area is 8.\\

\noindent Step 9. The square annulus $D_{9}$ comes from the edge $\overline{v_{18}v_{19}}$.  Its largest side length is three and the shortest side length is 1 (Figure \ref{piece} part $i$). Hence, its area is 8.\\

\noindent Step 10. The square cylinder $D_{10}$ is obtained from  $\overline{v_{19}v_{20}}$ and  $\overline{v_{20}v_{21}}$. Its base length is three, with a height of two (Figure \ref{piece} part $j$). Therefore, its area is 24.\\

\noindent Step 11. The square$D_{11}$ is constructed via the edges $\overline{v_{21}v_{22}}$ and $\overline{v_{22}v_{23}}$.  Its length side is three (Figure \ref{piece} part $k$). So, its area is 9.\\

\noindent We would like to point out that the pieces in Figure \ref{piece} parts $a$, $b$,  $c$, $e$ and $g$ are drawn with respect to  $x$, $z$, $t$-coordinates; while the remaining pieces are drawn with respect to $y$, $z$, $t$-coordinates.\\ 

\begin{main}\label{main}
There is a weakly minimal cubical 2-knot with area  130 isotopic to the spun trefoil knot.\\
\end{main}

\noindent \textit{Proof}.  If we apply a $(M1)$-move to $RC\Spin({\cal{A'}})$, then its area will remain equal. Hence, to reduce its area, we must apply the inverse of an $(M1)$-move, but by the spin construction, we have that $RC\Spin({\cal{A'}})$ is obtained by gluing disks, annuli, and cylinders. Some cylinders have radius $\frac{1}{2}$ and are deformed into square cylinders of base side one, hence there does not exist an integer 
$k >1$ such that $RC\Spin({\cal{A'}})$ is contained in the scaffolding of the corresponding cubulation ${\cal{C}}^k=H_k({\cal{C}})$, where $H_k:\mathbb{R}^4\rightarrow\mathbb{R}^4$ is the homothetic transformation given by $H_k(x)=kx$. Therefore, we can not apply the inverse of an $(M1)$-move to reduce its area. Now we will consider applying a $(M2)$-move. \\

\begin{figure}
\centering
\includegraphics[height=2.5cm]{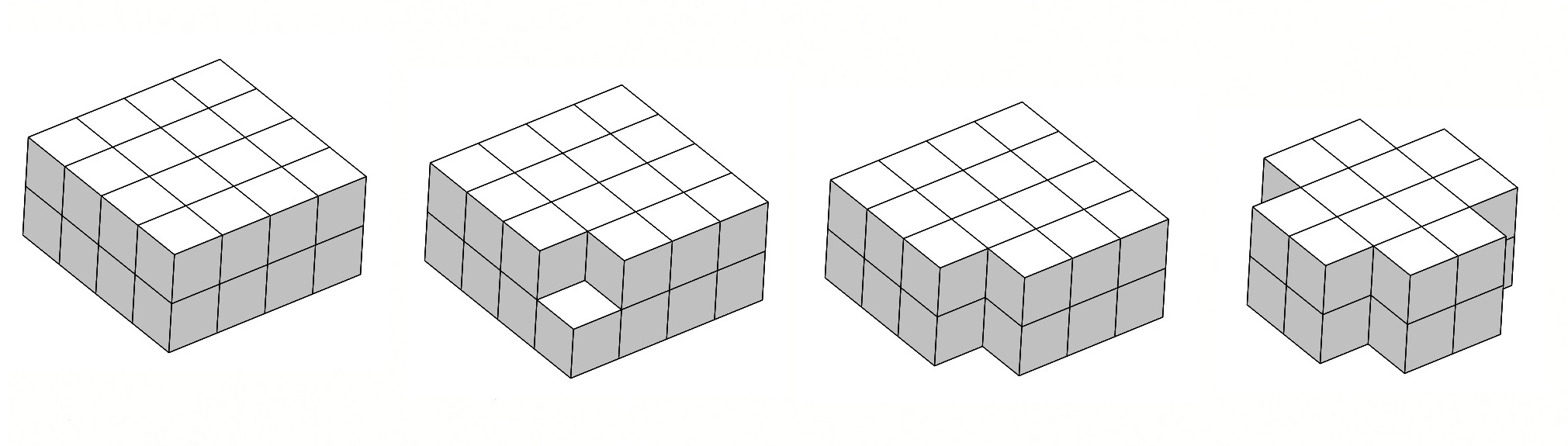}
\caption{A sequence of $(M2)$-moves applied to cubes sharing a corner-edge.}
\label{C2}
\end{figure}

\noindent If  $RC\Spin({\cal{A'}})$ intersects to a cube $Q$ of the canonical cubulation ${\cal{C}}$, then $RC\Spin({\cal{A'}})\cap Q$ can be either a vertex, an edge, one, two or three faces. However, to reduce the area of $RC\Spin({\cal{A'}})$, we must apply  $(M2)$-moves on cubes of ${\cal{C}}$ containing at least three faces of $RC\Spin({\cal{A'}})$. For instance, we can apply $(M2)$-moves on all the cubes intersecting a given edge $e\in RC\Spin({\cal{A'}})$, whose each of its end-points is a vertex of three faces of $RC\Spin({\cal{A'}})$ (corner vertices),  then after applying two $(M2)$-moves, we get a new cubical 2-knot $K^2$ such that $\Area (K^2)=\Area (RC\Spin({\cal{A'}}))-2$ and the edge $e$  does not belong to $K^2$, see Figure \ref{C2} and compare \cite{BH}. We will analyze these vertices for our 2-knot.\\

\noindent By construction, we have that $RC\Spin({\cal{A'}})=C\Spin(A')_{2}+b$, where $b=(0,0,\frac{1}{2},\frac{1}{2})$. Observe that both cubical 2-knots have the same area, since translation maps preserve areas, \emph{i.e.,} $\Area (RC\Spin({\cal{A'}}))=\Area (C\Spin(A')_{2})$. For simplicity, we will consider $C\Spin(A')_{2}$.
Notice that we have four cubical arcs $L_1, L_2, L_3$ and $L_4 \subset C\Spin(A')_{2}$ that are the images under the cubical process of
$R_{\theta_i}({\cal{A}})=\{(x_{1},x_{2},x_{3}\cos (2\pi \theta_i),x_{3}\sin(2\pi \theta_i))\,:\, (x_1,x_2,x_3)\in {\cal{A}}\}$,  for 
$\theta_i=0,\frac{1}{4},\frac{1}{2}, \frac{3}{4}$. In particular, $L_1={\cal{A'}}$ and each $L_i$ is contained into a copy of $\mathbb{R}^3_+$, say $\mathbb{R}^3_{+, i}$, which is a coordinate space, hence it has a cubulation ${\cal{C}}_{2,i}$ given by the restriction of the cubulation ${\cal{C}}_{2}$ of $\mathbb{R}^4$ to it (compare Theorem 1 in 
\cite{BH}).\\

\noindent \emph{Case 1.} The edge $e$ belongs to some arc $L_i$.  Since the construction of $C\Spin(A')_{2}$ is symmetric with respect to the corresponding coordinate spaces, we can assume that $L_i=L_1={\cal{A'}}$. Let $Q$ be a unit cube such that $e\subset Q$ and $Q+b$ is a cube of $\mathbb{C}$. Suppose that we have applied a $(M2)$-move to $C\Spin(A)_{2}$ to obtain a new cubical 2-knot $K^2$ such that $e$ does not belong to it, so the arc $L_1$ has been transformed into a new arc $L'_1\subset K^2$. This implies that if we remove a face containing $e$, then it must be replaced by at least one square, so by Theorem 3, the length of $L'_1$ is greater than or equal to the length of $L_1$; 
that is, $\Long (L_1')\geq \Long (L_1)$ hence $\Area (K^2)\geq \Area (RC\Spin({\cal{A'}}))$.\\

\noindent Suppose that $\Long (L_1')= \Long (L_1)$. This occurs if we apply a $(M2)$-move on a cube $Q$ containing two adjacent edges of $L_1$, then the remaining two edges on the same face of $Q$ replace these edges. Notice that these edges have length one, since we are considering $(M2)$-moves on faces of a cube $Q$ of side one, hence they correspond to $(M2)$-moves on $RC\Spin({\cal{A'}})=C\Spin(A')_{2}+b$. 
The above happens on the vertices $(1,-3,\frac{3}{2},0)$, $(-1,-3,\frac{3}{2},0)$, $(-1,-3,\frac{1}{2},0)$, $(-1,-1,\frac{1}{2},0)$,
$(2,-1,\frac{1}{2},0)$, $(2,-4,\frac{1}{2},0)$, $(0,-4,\frac{1}{2},0)$, $(0,-2,\frac{1}{2},0)$, $(0,-2,\frac{3}{2},0)$, and $(0,0,\frac{3}{2},0)$. However, by Theorem 3 we can not realise a $(M2)$- move on ${\cal{A}}$ to reduce the area of $RC\Spin({\cal{A'}})$. More specifically, 
If we realise a $(M2)$-move on the edges sharing the vertices $(1,-3,\frac{3}{2},0)$, $(-1,-3,\frac{3}{2},0)$, $(-1,-1,\frac{1}{2},0)$,
$(2,-4,\frac{1}{2},0)$, $(0,-4,\frac{1}{2},0)$, $(0,-2,\frac{1}{2},0)$, $(0,-2,\frac{3}{2},0)$, and $(0,0,\frac{3}{2},0)$, then the new arc $L'_1$ would have self-intersections. If 
 we apply an $(M2)$-move on the vertex $(-1,-3,\frac{1}{2},0)$, by Corollary 1, we have that $\Area (K^2) > \Area (RC\Spin({\cal{A'}}))$. Finally, on the vertex 
$(2,-1,\frac{1}{2},0)$ would have that $\Area (K^2)= \Area (RC\Spin({\cal{A'}}))$. \\
 
\noindent Summarizing, we obtain that $\Area (K^2)\geq \Area (RC\Spin({\cal{A'}}))$.\\
 
\noindent \emph{Case 2.} The intersection $e\cap L_i$ is a vertex $v_i\in C\Spin(A')_{2}$. As above, if we replace the vertex $v_i$ by applying a $(M2)$-move on 
$C\Spin(A)_{2}$, we obtain a new arc $L'_i$ whose length is greater than or equal to the length of $L_i$, to reduce the area of $C\Spin(A')_{2}$, both lengths must be equal. By the discussion on the previous case, we have only two possibilities: one is to replace the vertex $(-1,-3,\frac{1}{2},0)\in L_1$ (or the corresponding one on $L_i$) by the vertex $(-1,-2,\frac{3}{2},0)$ and the other is to replace the vertex  $(2,-1,\frac{1}{2},0)$ by the vertex  $(1,-2,\frac{1}{2},0)$. In both cases, by Theorem 3, we get a cubical 2-knot $K^2$ such that $\Area (K^2)\geq \Area (RC\Spin({\cal{A'}}))$.\\

\noindent \emph{Case 3.} The intersection $e\cap L_i$ is empty. Since $C\Spin(A')_{2}$ is obtained by gluing squares, square annuli, and square cylinders, and the side lengths of the squares can be only one or three, the smallest and biggest lengths of a square annulus are 1 and 3, respectively. The length of the base of a square cylinder can be only one or three. Suppose that $e$ is parallel to an edge $e_i$ of $L_1$ and both belong to the same face $Q$. If  we apply a 
 $(M2)$-move to remove $e_i$, then we obtain a new knot $K^2$ such that instead of $Q$ we get a new face $Q'$ parallel to it, hence our arc $L_1$ is replace by a new arc $L_1'\subset K^2$ and by the discussion on Case 1, we have that $\Area (K^2)= \Area (RC\Spin({\cal{A'}}))$. The same argument holds if $e_i$ and an end-point of $e$ belong to the same face $Q$.\\
 
\noindent Suppose that  $Q$ is a face such that $L_i\cap Q\neq\emptyset$, so $Q\cap e=\emptyset$. Then, the edge $e$ must belong to an arc $N_i$ that is the image under the cubical process of $R_{\theta_i}({\cal{A'}})=\{(x_{1},x_{2},x_{3}\cos (2\pi \theta_i),x_{3}\sin(2\pi \theta_i))\,:\, (x_1,x_2,x_3)\in {\cal{A'}}\}+b$,  for 
$\theta_i=\frac{1}{8},\frac{3}{8},\frac{5}{8}, \frac{7}{8}$ and $i=1,2,3,4$. Observe that $e+b$ is an edge of an arc $N_i$ for some $i$, and $N_i$ is contained in the canonical cubulation ${\cal{C}}$ of $\mathbb{R}^4$, more specific it is a copy of our arc ${\cal{A}}$ contained into a copy of $\mathbb{R}^3_+$, $R_{\theta_i}(\mathbb{R}^3_+)+n$ (see \cite{zeeman1}). By the same argument of Case 1, we have that if we apply a $(M2)$-move on $RC\Spin({\cal{A'}})=C\Spin(A')_{2}+b$ to obtain a new cubical 2-knot $K^2$ such that $e\not\subset K^2$, then $N_i$ is replaced by a new arc $N'_i$ satisfying that $\Long (N_i')\geq \Long (N_i)$ and as a consequence 
$\Area (K^2)\geq \Area (RC\Spin({\cal{A'}}))$.\\

\noindent  Consider the case  $\Long (N_i')=\Long (N_i)$. If we apply a $(M_2)$-move on a cube $Q$ containing two adjancent edges of $N_i$, then the remaining two edges are on the same face of $Q$, this implies that $N'_i$ is contained in the same copy of $\mathbb{R}^3_+$, so by the arguments of Case 1, we can conclude that there is no 
$(M2)$-move that reduces the area of $RC\Spin({\cal{A'}})=C\Spin(A')_{2}+b$. Therefore, the result follows. $\square$\\

\begin{main4}
The minimal area of the spun trefoil knot is smaller than or equal to 130; in other words,
$$
A (\Spin (T_{2,3}))\leq 130.
$$
\end{main4}


\vskip .5cm
\noindent{\large \bf Funding}  Ana Baray  and Juan Jos\'e Catal\'an were partially supported by SECIHTI-CONAHCyT (M\'exico) PhD Scholarship.\\

\noindent A. Baray. {\tt Centro de Investigaci\'on en Ciencias}. Instituto de Investigaci\'on en Ciencias B\'asicas y Aplicadas. Universidad Aut\'onoma del Estado de Morelos. Av. Universidad 1001, Col. Chamilpa.
Cuernavaca, Morelos, M\'exico, 62209. 

\noindent {\it E-mail address:} barayana@hotmail.com

\vskip .3cm

\noindent J. J. Catal\'an. {\tt Centro de Investigaci\'on en Ciencias}. Instituto de Investigaci\'on en Ciencias B\'asicas y Aplicadas. Universidad Aut\'onoma del Estado de Morelos. Av. Universidad 1001, Col. Chamilpa.
Cuernavaca, Morelos, M\'exico, 62209. 

\noindent {\it E-mail address:} josecatalan103@gmail.com

\vskip .3cm

\noindent G. Hinojosa. {\tt Centro de Investigaci\'on en Ciencias}. Instituto de Investigaci\'on en Ciencias B\'asicas y Aplicadas. Universidad Aut\'onoma del Estado de Morelos. Av. Universidad 1001, Col. Chamilpa.
Cuernavaca, Morelos, M\'exico, 62209. 

\noindent {\it E-mail address:} gabriela@uaem.mx 

\vskip .3cm

\noindent R. Valdez. {\tt Centro de Investigaci\'on en Ciencias}. Instituto de Investigaci\'on en Ciencias B\'asicas y Aplicadas. Universidad Aut\'onoma del Estado de Morelos. Av. Universidad 1001, Col. Chamilpa.
Cuernavaca, Morelos, M\'exico, 62209. 

\noindent {\it E-mail address:} valdez@uaem.mx 

\end{document}